\documentclass[a4paper,10pt]{article}
\usepackage[utf8x]{inputenc}
\usepackage{graphicx}

\usepackage{amssymb}
\usepackage{amsmath}
\usepackage{amsfonts}
\usepackage{multirow}
\usepackage{moreverb}
\usepackage{subfig}
\usepackage{rotating}

\usepackage[colorlinks,bookmarksopen,bookmarksnumbered,citecolor=red,urlcolor=red]{hyperref}
\hypersetup{pdftitle={A "well-balanced" finite volume scheme for blood flow simulation},
pdftoolbar=true,
pdfmenubar=true,
pdfauthor={O. Delestre, P.-Y. Lagree},
pdfsubject={applied mathematics},
pdfcreator={Delestre Olivier},
pdfproducer={Delestre Olivier},
pdfkeywords={blood flow simulation} {well-balanced scheme} {finite volume scheme} {hydrostatic reconstruction}
 {man at eternal rest} {semi-analytical solutions} {shallow water} {Saint-Venant} {friction}}

\title{A "well-balanced" finite volume scheme for blood flow simulation}
\author{{{O. Delestre}\footnote{Corresponding author: Delestre@unice.fr, presently at: Laboratoire de Math\'ematiques
 J.A. Dieudonn\'e  (UMR CNRS 7351) -- Polytech Nice-Sophia , Universit\'e de Nice -- Sophia Antipolis, Parc Valrose, 06108 Nice cedex 02,
 France}$^{\; ,}$}\footnote{CNRS \& UPMC Universit\'e Paris 06, UMR 7190, 4 place Jussieu,
Institut Jean Le Rond d'Alembert, Bo\^ite 162, F-75005 Paris, France}, P.-Y. Lagr\'ee$^{\dag}$}

\begin{document}
\maketitle

\begin{abstract}
We are interested in simulating blood flow in arteries  with a one dimensional model. Thanks to recent developments in the
 analysis of hyperbolic system of conservation laws (in the Saint-Venant/ shallow water equations context) we will perform a
 simple finite volume scheme. We focus on conservation properties of this scheme which were not previously considered. To emphasize
 the necessity of this scheme, we present how a too simple numerical scheme may induce spurious flows when the basic static shape
 of the radius changes. On contrary, the proposed scheme is "well-balanced":  it preserves equilibria of $Q=0$. Then examples of
 analytical or linearized solutions with and without viscous damping are presented to validate the calculations. 
 The influence of abrupt change of basic radius is emphasized in the case of an aneurism.
\end{abstract}


\section{Introduction}
\vspace{-2pt}

As quoted by Xiu and Sherwin  \cite{Xiu07} the one dimensional system of equations for blood flow in arteries was written 
long time ago by Leonhard Euler in 1775.  Of course, he noticed (see Parker \cite{Parker09}) that it was far too difficult to
 solve. 
Since then, the pulsatile wave flow has been understood and is explained in classical books such as Lighthill \cite{Lighthill78} and Pedley \cite{Pedley80}. 
More recently a  huge lot of work has been done and a lot of methods of resolution have been proposed to handle this set of
 equations, among them
\cite{VanSteenhoven86}, \cite{Zagzoule86}, \cite{Stergiopulos92},
 \cite{Olufsen99}, \cite{Lagree96}
 \cite{Lagree00} \cite{Sherwin03}, \cite{Xiu07}, \cite{Fernandez04b}, \cite{Wibmer04}, \cite{Martin05}, \cite{Formaggia06},
 \cite{Cavallini08},
 \cite{Alastruey08}, \cite{Cavallini10}, \cite{Willemet11}, \cite{Willemet11,Marchandise08}... 
 and we forgot a lot of actors.  
 Even if recent progress of fluid structure interaction in 3D have been performed (Gerbeau {\it et. al.} \cite{Fernandez07},
 Van de Vosse {\it et. al.}  \cite{vandeVosse03b},  among others),
 those advances have not made the 1-D modelisation obsolete. On the contrary, it is needed for validation, for physical
 comprehension, for boundary conditions, fast or real time computations and it is always relevant for full complex networks of
 arteries (\cite{Zagzoule86}, \cite{Stergiopulos92},  \cite{Saito11}, \cite{Olufsen99,Olufsen00}) or veins (\cite{Fullana09}) or in micro
 circulation in the brain \cite{Pindera09}.

In most of the previous references it is solved thanks to the  method of characteristics,  to finite differences or finite
 volumes (that we will discuss and improve in this paper), or finite elements methods (Galerkin discontinuous).
 Each method has its own advantages.

In the meanwhile, the set of shallow water equations (though younger, as the equations have been written  by Adh\'emar Barr\'e
 de Saint-Venant in 1871 \cite{saintvenant71}) have received a large audience as well. They are used for modelling the flow in
 rivers \cite{Goutal02}, \cite{Burguete08} (and networks of rivers), in lakes, rain overland flow \cite{Delestre09b,Delestre09},
 \cite{Tatard08}, in dams \cite{Valiani99a,Valiani02} or the long waves in shallow water seas (tides in the Channel for example
 \cite{Sampson05}, or for Tsunami modelling \cite{Dutykh07}, \cite{Popinet11}). 
It has been observed that Saint-Venant equations with source terms (which are the shape of the topography, and the viscous
 terms) present numerical difficulties for steady states. The splitting introduced by the discretisation creates an extra
 unphysical flow driven by the change of the topography.
A configuration with no flow, for example a lake, will not stay at rest: thus, the so called equilibrium of the "lake at rest"
 is not preserved. So new schemes, called "well-balanced" have been constructed to preserve this equilibrium. Bermudez and
 Vazquez \cite{Bermudez94,Bermudez98} have modified the Roe solver to preserve steady states. Other ways to adapt exact  or approximate
 Riemann solver to non-homogeneous case have been proposed by LeVeque \cite{LeVeque98} and Jin \cite{Jin01a}.
Following the idea of the pioneer work of Greenberg and LeRoux \cite{Greenberg96}, \cite{Gosse00} and \cite{Gallouet03} have
 proposed schemes based on the solution of the Riemann problem associated with a larger system (a third equation on the bottom
 topography variable is added). A central scheme approach has been developped in \cite{Kurganov02}. In Perthame and Simeoni
 \cite{Perthame01}, the source terms are included in the kinetic formulation. In \cite{Audusse04c} and \cite{Bouchut04},
 a well-balanced hydrostatic reconstruction has been derived. It can be applied to any conservative finite volume scheme
 approximating the homogeneous Saint-Venant system. In \cite{Vazquez-Cendon99}, \cite{Chinnayya99}, \cite{Bouchut04} and
 \cite{LukacovaMedvidova05b,LukacovaMedvidova06} the friction source term is treated under a well-balanced way. More recently,
 in \cite{Lee10}, the system of shallow water with topography is rewritten under a homogeneous form. This enables them to
 get Lax-Wendroff and MacCormack well-balanced schemes. With well-balanced methods, spurious current does not appear.
 
Thus we will follow those authors and construct by analogy a similar method but for flows in elastic tubes rather than
 for free surface flows. Even if for artery flows the wave behaviour of the solution is very important (forced regime with
 lot of reflexions), we have in mind complex networks. In this case we will have to deal from the arteries to the veins with
 eventually microcirculation in the brain. In the smaller and softer vessels, the 1D approximation is still valid, the
 pulsatile behaviour is less important, and there are lot of changes of sections (across stenosis, aneurisms, valves).
 Those changes of sections have to be treated numerically with great care, which has not been the case up to now in the
 litterature.

 So here, we will follow the authors of the Saint-Venant community and construct by analogy a similar method but for artery
 flow. At first we derive the equations
 written in conservative form, looking in the literature, we found that it has been only written by \cite{Wibmer04} in a PhD
 thesis and by \cite{Cavallini08}, \cite{Cavallini10} but not exploited in the "well-balanced" point of view.
Mostly (\cite{Xiu07}, \cite{Marchandise08}), equations are in fact written in a non-conservative form. 
With non-conservative schemes, we have consistency defaults and shocks are not uniquely defined, see among others \cite{Hou94},
 \cite{Berthon07} and \cite{Castro08b}. Furthermore, the discharge $Q$ is not conserved. Not a lot of attention has been
 attached to this problem up to now. But this non conservation is a real problem if we consider networks with tubes with
 non uniform section and ultimately quasi steady flow in small terminal vessels.

Hence after the introduction of the set of equations in section \ref{model} and deriving with this new point of view
the conservative equations, we construct by analogy the numerical scheme both at first and second order in section \ref{method}. Then section
 \ref{validation} is devoted to several test cases to validate the different parts of the equations. The most important case
 is the "man at eternal rest" by analogy to the "lake at rest" which would not be preserved by non-conservative schemes.
The aim of this paper is to use this property to construct efficient methods of resolution.


\section{The 1D model}\label{model}
\subsection{Derivation of the equations}
We first derive the model equations, they are of course now classical, as already mentioned. 
Starting from unsteady incompressible axi-symmetrical Navier-Stokes equations, and doing a long wave approximation one
 obtains a first  set of Reduced Navier-Stokes equations (\cite{Lagree00}, \cite{Lagree05}). 
One has the longitudinal convective term, the longitudinal pressure gradient, and only the transverse viscous term, the
 longitudinal being negligible.
The pressure is constant across the cross section and depends only of $x$ the longitudinal variable.
The incompressibility is preserved.

\begin{center}
\begin{figure}[htbp]
   \includegraphics[width=12.3cm]{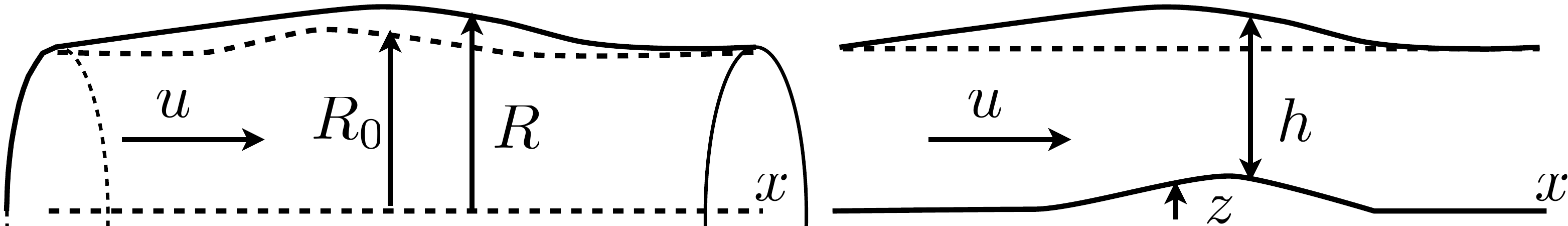} 
      \caption{1D models: left the  blood flow model ($u$ is the mean velocity, $R$ the radius, $A=\pi R^2$ the
 instantaneous area, $R_0$ the radius at rest, $A_0=\pi R_0^2$ the area with no flow), right the shallow water model
      ($u$ is the mean velocity, $h$ the water column height, $z$ the topography).}
   \label{fig:s1}
\end{figure}
\end{center}

Then, integrating over the section, one obtains two integrated equations (kind of Von K\'arm\'an equation). The first one is
 obtained from incompressibility and application of boundary condition and relates the variation of section to the variations
 of flow. The second one is obtained from the momentum. 
The final equations are not closed and one has to add some hypothesis on the shape of the velocity profile to obtain a closed
 system (some discussions of this are in \cite{Lagree00}).  We then obtain  the following system with dimensions 
which is the 1-D model of flow:
\begin{equation}
 \left\{\begin{array}{l}
         \partial_t A+\partial_x Q=0\\
	\partial_t Q+\partial_x \dfrac{Q^2}{A}+\dfrac{A}{\rho}\partial_x p=-C_f\dfrac{Q}{A}
        \end{array}\right.,\label{system1}
\end{equation}
where
\(A(x,t)\) is the cross-section area ($A=\pi R^2$, where $R$ is the radius of the vessel), \(Q(x,t)=A(x,t)u(x,t)\) is the flow
 rate or the discharge, \(u(x,t)\) the mean flow velocity, and  \(\rho\) the blood density (see figure \ref{fig:s1}).
The first equation is without any approximation.
An extra factor of value $4/3$ may appear in front of $Q^2/A$ if the chosen profile is a Poiseuille profile (with
either large viscosity $\nu$ or low frequency $\omega/(2 \pi)$ so that   $\alpha= R_0\sqrt{\omega/\nu}$ the  Womersley number
 is small).
This comes from the fact that this term is in fact an approximation of $\int_0^R2 \pi u(t,x,r)^2 r dr$ (this $u(t,x,r)$ being
 the actual longitudinal velocity, not to be confused with $u(x,t)=\int_0^R2 \pi u(t,x,r) r dr$). The chosen unit value
 corresponds to a flat profile. The viscous term has been modelised by a friction proportional to $Q/A$ (see
 \cite{Zagzoule91} \cite{Lagree00}):
\(C_f=8\pi\nu\) where $\nu$ is the blood viscosity.  
Again, a friction term written 
$-8\pi\nu Q/A$  is an approximation of the effective term which is $\nu R \partial u/\partial r |_R$ (with here again
  $u(t,x,r)$ the real longitudinal velocity). This specific choice corresponds to a Poiseuille flow (small $\alpha$), and is
 therefore non consistent with the coefficient in the convective term. This is not a real problem, those coefficients are
 adjusted in the literature. 
Note that in Saint-Venant community, the skin friction modelling involves a square of $Q$, reminiscent to the turbulent nature
 of river flow. This is seldom used in blood-flow.
To conclude on the mechanical modelling we will take the form of Eq. \ref{system1} for granted and we will no more discuss it.

Finally, one has to model the arterial wall, the pressure has to be linked with the displacement of the vessel. A simple elastic
 linear term in the variation of radius is taken:
\begin{equation}
 p-p_0=k (R-R_0) \mbox{  or as well } p=p_0+k\dfrac{\sqrt{A}-\sqrt{A_0}}{\sqrt{\pi}}\label{pressure},
\end{equation}
 \(A_0(x)\) is the
 cross section at rest ($A_0=\pi R_0^2$, where $R_0$ is the radius of the vessel which may be variable in the case of aneurism,
 stenosis or taper), and the stiffness  \(k\) is taken as a constant.  The external pressure $p_0$ is supposed constant as
 well. Again, this relation which is here a simple elastic string model may be enhanced. Non linear terms,
 non local (tension) terms and/or dissipative unsteady terms may be introduced (viscoelasticity). We take here
 the most simple law, but we are aware of the  possibilities of complexity. 
 
Note that the derivation of those equations may be done  with only 1D fluid mechanics arguments (\cite{Pindera09} or
 \cite{Kundu04}). One has to be very careful then to write precisely the action of the pressure on the lateral walls (for
 example Kundu  \cite{Kundu08} page 793 does a mistake in his text book). But, the real point of view to derive the equations
 is the one explained at the beginning of this section and involves the Reduced Navier-Stokes equations (\cite{Lagree00}).

We will now change a bit this form to write a real conservative form. Note that authors (\cite{Sherwin03},
 \cite{Marchandise08}) write system  \eqref{system1} as 
$$\partial_tU+\partial_x F(U)= S(U),$$
 they use the variables $U=(A,u)$, so  that they identify a kind of flux $F(U) =(Au,u^2/2+p/\rho)$ and a source term
 $S(U)=  -C_f{u}/{A}$.
But $u$ is not a conservative quantity, $Q$ is, as we will see next.

Of course, those Saint-Venant and blood flow  equations are very similar to the 1D Euler compressible equations that one can
 write in nozzles, or in acoustics (Lighthill \cite{Lighthill78}), in this case the density $\rho$ is variable and linked to
 the pressure by the isentropic law  $p \propto \rho^\gamma$. This relation is the counter part of the relation of pressure
 which are in Saint-Venant and respectively in the artery:
hydrostatic balance $p=\rho g \eta$ (density, gravitational acceleration, level of the free surface) and respectively
  the elasticity response of the artery $p= k (R-R_0)$ (stiffness of the wall, current radius of the vessel, reference radius
 at rest).   
  
\subsection{Conservative system}
So, to write the conservative system,  the second equation of \eqref{system1} is rewritten and developed. 
Thus we get the following system of conservation laws:
\begin{equation}
 \left\{\begin{array}{l}
         \partial_t A+\partial_x Q=0\\
	\partial_t Q+\partial_x \left[\dfrac{Q^2}{A}+\dfrac{k}{3\rho\sqrt{\pi}}A^{3/2}\right]
=\dfrac{kA}{2\rho\sqrt{\pi}\sqrt{A_0}}\partial_x A_0-C_f\dfrac{Q}{A}
        \end{array}\right.,\label{system2}
\end{equation}
we can write \eqref{system2} under a more compact form:
\begin{equation}
 \partial_t U+\partial_x F(U)=S(U), \label{compact-form}
\end{equation}
with a new definition of $U$, $F(U)$ and $S(U)$.  We identify  \(U\) the vector of the conservative variables, $F(U)$ the flux
\begin{equation}
 U=\left(\begin{array}{c}
          A\\
	  Q
         \end{array}\right),\;
F(U)=\left(\begin{array}{c}
            Q\\
	    \dfrac{Q^2}{A}+\dfrac{k}{3\rho\sqrt{\pi}}A^{3/2}
           \end{array}\right)\label{compact-form2}
\end{equation}
and the source term takes into account the initial  shape of the vessel $A_0$ and the friction term:
\begin{equation}
\begin{array}{cc}
 S(U)     & =\left(\begin{array}{c}
             0\\
	     \dfrac{kA}{\rho\sqrt{\pi}}\partial_x \sqrt{A_0}-C_f\dfrac{Q}{A}
            \end{array}\right).
\end{array}
\label{compact-form3}
\end{equation}
This system is strictly hyperbolic when $A>0$ (which should be the case in arteries). Indeed, this means that we have
\begin{equation*}
 \partial_x F(U)=\left(\begin{array}{cc}
                        0 							& 1\\
			\dfrac{k\sqrt{A}}{2\rho\sqrt{\pi}}-\dfrac{Q^2}{A^2}	& \dfrac{2Q}{A}
                       \end{array}\right).
\partial_x \left(\begin{array}{c}
       A\\
       Q
      \end{array}\right)
=J(U).\partial_x U,
\end{equation*}
where the Jacobian matrix \(J(U)\) has two real eigenvalues (which is the definition of hyperbolicity)
\begin{equation}
\lambda_1 =\dfrac{Q}{A}-\sqrt{\dfrac{k\sqrt{A}}{2\rho\sqrt{\pi}}}=u-c
  \quad\text{and}\quad
 \lambda_2 =\dfrac{Q}{A}+\sqrt{\dfrac{k\sqrt{A}}{2\rho\sqrt{\pi}}}=u+c.
\end{equation}
We recognise \(c\) as the well-known Moens Korteweg wave propagation speed (\cite{Pedley80}):
\begin{equation*}
 c=\sqrt{\dfrac{k\sqrt{A}}{2\rho \sqrt{\pi}}}.
\end{equation*}
The system \eqref{system2} without friction {\it i.e.} with $C_f=0$ 
admits special solutions which are steady state solutions, that means that we have
\begin{equation}
 \partial_t A=\partial_t u=\partial_t Q=0.
\end{equation}
Thus the mass conservation equation gives the discharge conservation:  $\partial_xQ(x,t)=0$. 
The momentum conservation equation, under the steady state and frictionless hypothesis
reduces to 
\begin{equation}
 \partial_x\left[\dfrac{{Q_0}^2}{2A^2}+b\sqrt{A}-b\sqrt{A_0}\right]=0,
\end{equation}
with $b=k/(\rho\sqrt{\pi})$. So that system \eqref{system2}  integrates to the two constants:
\begin{equation}
\left\{\begin{array}{l}
 Q(x,t)=Q_0\\
 \dfrac{{Q_0}^2}{2A^2}+b\sqrt{A}-b\sqrt{A_0}=Cst 
       \end{array}\right..\label{bernoulli-BF}
\end{equation}
The first is a constant flux, the second is the Bernoulli constant.

We have of course  made analogy between \eqref{system2} and the one-dimensional Saint-Venant (or shallow water) system:
\begin{equation}
\left\{\begin{array}{l}
 \partial_t h+\partial_x q =0\\
\partial_t q + \partial_x \left(\dfrac{q^2}{h}+\dfrac{gh^2}{2}\right)=-gh\partial_x z -c |q|\dfrac{q}{h^\beta}
\end{array}\right.,\label{SW}
\end{equation}
where $h(x,t)$ is the water height, $u(x,t)$ the mean flow velocity, $q(x,t)=h(x,t)u(x,t)$ the flow rate and $z(x,t)$ the
 topography of the bottom and $c$ is the friction coefficient. With $\beta=2$ we recover Chezy's or Darcy-Weisbach's law
 (depending on the way $c$ is written) and with $\beta=7/3$, we get the Manning-Strickler's friction law.
 
This system admits steady state frictionless solutions too:
\begin{equation}
 \partial_t h=\partial_t u=\partial_t q=0.
\end{equation}
Thus we get the conservation of discharge and the Bernoulli's relation
\begin{equation}
 \left\{\begin{array}{l}
	  q=q_0\\
	  \dfrac{{q_0}^2}{2gh^2}+z+h=Cst
        \end{array}\right.,\label{bernoulli-SV}
\end{equation}
which is the exact analogous of \eqref{bernoulli-BF}. In the literature we can find numerical schemes preserving steady states
 under the form \eqref{bernoulli-SV}, but they are complex to handle (see \cite{Castro07}, \cite{Noelle07}, \cite{Thanh08}
 and \cite{Bouchut09}). Most of the time, more simple steady states are
 considered
\begin{equation}
 \left\{\begin{array}{l}
         q=u=0\\
	\partial_x \left(g\dfrac{h^2}{2}\right)+gh\partial_x z=0
        \end{array}\right.,\label{steady-hydro}
\end{equation}
which correspond from a mechanical point of view to the modelling of a "lake at rest". We have a balance called hydrostatic
 balance between the hydrostatic pressure and the
 gravitational acceleration down an inclined bottom $z$. The second term of \eqref{steady-hydro} reduces to
\begin{equation}
 H=h+z=Cst,
\end{equation}
where $H$ is the water level. Since the work of \cite{Greenberg96}, schemes preserving at least \eqref{steady-hydro} are called
 well-balanced schemes. This allows to have an efficient treatment of the source term. Thus, in practice, schemes preserving
 \eqref{steady-hydro} give good results even in unsteady  cases. 
 The analogous of the "lake at  rest" for the system \eqref{system2} without friction,
can be called the "man at eternal rest" or "dead man equilibrium", it   writes
\begin{equation}
 \left\{\begin{array}{l}
         Q=u=0\\
	\partial_x \left(\dfrac{k}{3\rho \sqrt{\pi}}A^{3/2}\right)-\dfrac{kA}{\rho \sqrt{\pi}}\partial_x \sqrt{A_0}=0
        \end{array}\right..\label{rest-bf}
\end{equation}
In this case we have $\sqrt{A}=\sqrt{A_0}$.
We will now use this property to construct the schemes.

\begin{figure}[htbp]
   \centering
   \includegraphics[width=5cm]{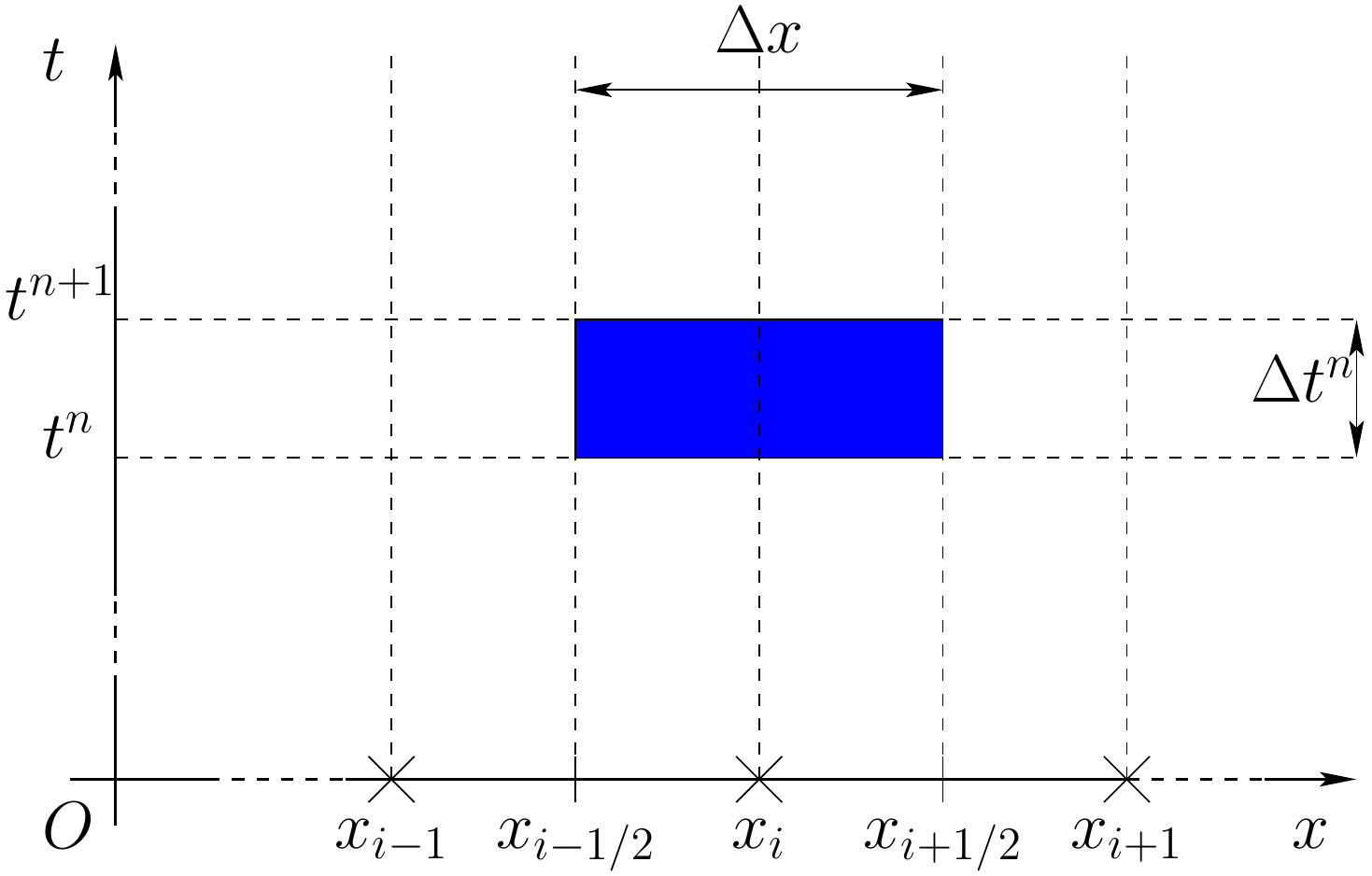} 
      \caption{The time and space stencil.}
   \label{fig:s2}
\end{figure}

\section{The numerical method} \label{method}
In this section, we describe the scheme for system \eqref{system1} at first and second order, with the space and time
 discretisation illustrated in Figure \ref{fig:s2}.
\subsection{Convective step}
For the homogeneous system {\it i.e.} system \eqref{compact-form} without source term: \(S(U)=0\), a first-order conservative
 finite volume scheme writes simply (see figure \ref{fig:s2}):
\begin{equation}
 \dfrac{U_i^{n+1}-U_i^{n}}{\Delta t}+\dfrac{F^n_{i+1/2}-F^n_{i-1/2}}{\Delta x}=0, \label{scheme}
\end{equation}
where \(U_i^n\) is an approximation of \(U\)
\begin{equation*}
 U_i^n\simeq \dfrac{1}{\Delta x} \int_{x_{i-1/2}}^{x_{i+1/2}} U(x,t_n)dx
\end{equation*}
 \(n\) refers to time \(t_n\) with \(t_{n+1}-t_{n}=\Delta t\) and \(i\) to the cell
 \(C_i=(x_{i-1/2},x_{i+1/2})=(x_{i-1/2},x_{i-1/2}+\Delta x)\). The two points numerical flux \(F_{i+1/2}\) is an
 approximation of the flux function \eqref{compact-form2} at the cell interface \(i+1/2\)
\begin{equation}
 F_{i+1/2}={\cal F}(U_i,U_{i+1}).
\end{equation}
This flux function will be detailed in section \ref{numerical-flux}.

\subsection{Source term treatment}

\subsubsection{"Topography" treatment}
\label{naive}
The source term $\dfrac{kA}{\rho \sqrt{\pi}}\partial_x \sqrt{A_0}$ is the analogous of the topography source term 
 ($-gh\partial_x z$) in the shallow water system \eqref{SW}.
 Looking at \eqref{system2}, it might be treated explicitly which writes
\begin{equation}
 U_i^{n+1}= U_i^{n}- \dfrac{\Delta t}{\Delta x}\left({F^n_{i+1/2}-F^n_{i-1/2}}\right)+\Delta t S(U_i^n) \label{scheme2}
\end{equation}
where the source is written simply by evaluating the derivative of $A_0$: 
\begin{equation*}
S(U_i^n)=\left(\begin{array}{c}
                0 \\
\dfrac{kA_i^n} {2\rho\sqrt{\pi}\sqrt{{A_0}_i}} \dfrac{{A_0}_{i+1/2}-{A_0}_{i-1/2}}{\Delta x} 
               \end{array}\right),
\end{equation*}
with the definition
\begin{equation*}
{A_0}_{i+1/2}=\dfrac{{A_0}_{i+1}+{A_0}_{i}}{2}. 
\end{equation*}
We will illustrate the fault of this naive method in sections \ref{deadman} and \ref{deadwave}.
 We will prefer a well-balanced method inspired
 from the hydrostatic reconstruction (see \cite{Audusse04c} and \cite{Bouchut04}). As for the hydrostatic reconstruction we do
 not consider the friction term and we are looking for a scheme which preserves steady states at rest \eqref{rest-bf},
 the second equation in \eqref{rest-bf} means that we have
\begin{equation}
 \sqrt{A}-\sqrt{A_0}=Cst,\label{steady-state2}
\end{equation}
as for the hydrostatic reconstruction we will use locally this relation to perform a reconstruction of the variable \(A\).
On each part of the interface, we have the following reconstructed variables
\begin{equation}
 \left\{\begin{array}{l}
         \sqrt{A_{i+1/2L}}=\max(\sqrt{A_i}+\min(\Delta\sqrt{A_0}_{i+1/2},0),0) \\
	U_{i+1/2L}=(A_{i+1/2L},A_{i+1/2L}u_i)^t \\
	 \sqrt{A_{i+1/2R}}=\max(\sqrt{A_{i+1}}-\max(\Delta\sqrt{A_0}_{i+1/2},0),0) \\
	U_{i+1/2R}=(A_{i+1/2R},A_{i+1/2R}u_{i+1})^t
        \end{array}\right.,\label{hydro-rec}
\end{equation}
where \(\Delta \sqrt{A_0}_{i+1/2}=\sqrt{{A_0}_{i+1}}-\sqrt{{A_0}_{i}}\).

For consistency the scheme \eqref{scheme} has to be slightly modified under the form
\begin{equation}
 U_i^{n+1}=U_i^{n}-\dfrac{\Delta t}{\Delta x}(F^n_{i+1/2L}-F^n_{i-1/2R}), \label{scheme-hydro}
\end{equation}
where
\begin{equation}
 \begin{array}{l}
    F^n_{i+1/2L}=F^n_{i+1/2}+S_{i+1/2L}\\
    F^n_{i-1/2R}=F^n_{i-1/2}+S_{i-1/2R}
 \end{array},
\end{equation}
with
\begin{equation}
\begin{array}{l}
 S_{i+1/2L}=\left(\begin{array}{c}
        0 \\
      {\cal P}(A_i^n)-{\cal P}(A_{i+1/2L}^n)
       \end{array}\right) \\
S_{i-1/2R}=\left(\begin{array}{c}
        0 \\
      {\cal P}(A_i^n)-{\cal P}(A_{i-1/2R}^n)
       \end{array}\right)
\end{array}
\end{equation}
and
\begin{equation}
 {\cal P}(A)=\dfrac{k}{3\rho\sqrt{\pi}}A^{3/2}. 
\end{equation}
The numerical flux $F^n_{i+1/2}$ is now calculated with the reconstructed variables \eqref{hydro-rec}:
\begin{equation*}
 F^n_{i+1/2}={\cal F}(U^n_{i+1/2L},U^n_{i+1/2R}).\label{num-flux}
\end{equation*}

Having treated the term linked to the variation of radius we now turn to the friction term.

\subsubsection{Friction treatment}

Concerning the friction treatment two methods are presented here, namely: a semi-implicit treatment (\cite{Fiedler00},
 \cite{Bristeau01}, and \cite{Liang09b}) which is classical in shallow water simulations and the apparent topography method
 (see \cite{Bouchut04}, \cite{MangeneyCastelnau05,Mangeney07}) which is a well-balanced method.

\paragraph{Semi-implicit treatment}

We use \eqref{scheme-hydro} as a prediction without friction, {\it i.e.}:
\begin{equation}
 U_i^{*}=U_i^{n}-\dfrac{\Delta t}{\Delta x}(F^n_{i+1/2L}-F^n_{i-1/2R})
\end{equation}
and a friction correction (via a semi-implicit treatment) is applied on the predicted variables ($U_i^{*}$):
\begin{equation}
 A_i^*\left(\dfrac{u_i^{n+1}-u_i^*}{\Delta t}\right)=-C_f u_i^{n+1},\label{semi-implicit}
\end{equation}
thus we get $u_i^{n+1}$ and we have $A_i^{n+1}=A_i^*$.
It preserves zero velocity.

\paragraph{Apparent topography}
The apparent topography method (see \cite{Bouchut04}) consists in writing the second
 equation of the system under the following form
\begin{equation}
 \partial_t Q+\partial_x \left[\dfrac{Q^2}{A}+\dfrac{k}{3\rho \sqrt{\pi}}A^{3/2}\right]
=\dfrac{kA}{\rho \sqrt{\pi}} \partial_x \sqrt{{A_0}_{app}},\label{system-app}
\end{equation}
with
\begin{equation*}
 \partial_x \sqrt{{A_0}_{app}}=\partial_x \left(\sqrt{A_0}+b\right),
\end{equation*}
where
\begin{equation*}
 \partial_x b=-\dfrac{\rho \sqrt{\pi} C_f}{k}\dfrac{Q}{A^2}.
\end{equation*}
The hydrostatic reconstruction \eqref{hydro-rec} is performed with the corrected term \(\sqrt{{A_0}_{app}}\).

\subsection{Numerical flux}\label{numerical-flux}
Several numerical fluxes might be used, some of them are defined in the following.
\paragraph{Rusanov flux}\label{rusanov-flux}
Following \cite{Bouchut04}, the Rusanov flux writes
\begin{equation*}{\cal F}(U_L,U_R)=
\dfrac{F(U_L)+F(U_R)}{2}-c\dfrac{U_R-U_L}{2},
\end{equation*}
with
\[c=\sup\limits_{U=U_L,U_R}({\sup\limits_{j\in\{1,2\}}} |\lambda_j(U)|),\]
where \(\lambda_1(U)\) and \(\lambda_2(U)\) are the eigenvalues of the system.

\paragraph{HLL flux}\label{HLL-flux}
Following \cite{Bouchut04}, the HLL flux (Harten, Lax and van Leer \cite{Harten83}) writes
\begin{equation}
 {\cal F}(U_L,U_R)=\left\{\begin{array}{ll}
F(U_L) & \text{if}\;0 \leq c_1 \\
\dfrac{c_2F(U_L)-c_1F(U_R)}{c_2-c_1}+\dfrac{c_1c_2}{c_2-c_1}(U_R-U_L) & \;\text{if}\;c_1<0<c_2 \\
F(U_R) & \text{if}\; c_2 \leq 0
\end{array}\right.,
\end{equation}
with
\[c_{1}={\inf\limits_{U=U_L,U_R}}({\inf\limits_{j\in\{1,2\}}}\lambda_{j}(U))\;
\text{and}\;c_{2}={\sup\limits_{U=U_L,U_R}}({\sup\limits_{j\in\{1,2\}}}\lambda_{j}(U)),\label{hll2}\]
where \(\lambda_1(U)\) and \(\lambda_2(U)\) are the eigenvalues of the system.

\paragraph{VFRoe-ncv flux}\label{VFRoe-flux}
Adapting the main line of \cite{Marche07b} and \cite{Marche08}, we get the following VFRoe-ncv flux for system \eqref{system2}
\begin{equation*}{\cal F}(U_L,U_R)=
\left\{\begin{array}{l}
F(U_L)\;\text{if}\; 0<\lambda_1(\widetilde U) \\
F(U_*)\;\text{if}\; \lambda_1(\widetilde U)<0<\lambda_2(\widetilde U) \\
F(U_R)\;\text{if}\; \lambda_2(\widetilde U)<0 \end{array}\right.,
\end{equation*}
with the mean state
\[\widetilde W=
\left(\begin{array}{c}
4\tilde c \\
\tilde u \end{array}\right)
=
\left(\begin{array}{c}
2\left(c_L+c_R\right) \\
\dfrac{u_L+u_R}{2} \end{array}\right)
\]
and
\begin{equation*}
 U_*=\left(\begin{array}{c}
            A_*\\
	    A_* u_*
           \end{array}\right),
\end{equation*}
where \(A_*\) et \(u_*\) are the Roe mean states, we get them thanks the Riemann invariant
\begin{equation}
 \left\{\begin{array}{l}
         u_R-4c_R=u_*-4c_*\\
	  u_L+4c_L=u_*+4c_*
        \end{array}\right.,
\end{equation}
thus we have
\begin{equation*}
\left\{\begin{array}{l}
u_*=\tilde u-2(c_R-c_L) \\
c_*=\tilde c-\dfrac{1}{8}(u_R-u_L)\end{array}\right..
\end{equation*}
When $\lambda_1(U_L)<0<\lambda_1(U_R)$ or $\lambda_2(U_L)<0<\lambda_2(U_R)$, we can get non entropic solutions. We perform
 an entropy correction thanks to the Rusanov flux.

\paragraph{Kinetic flux}\label{kinetic-flux}
Adapting the main line of \cite{Audusse99} and \cite{Audusse05}, we get the following kinetic flux for system \eqref{system2}
\begin{equation*}{\cal F}(U_L,U_R)=F_+(U_L)+F_-(U_R),
\end{equation*}
with
\begin{equation*}
 F_+(U_L)=\dfrac{A_L}{2\sqrt{3T_L}}
\left(\begin{array}{c}
       \dfrac{{M_+}^2-{M_-}^2}{2}\\
      \dfrac{{M_+}^3-{M_-}^3}{3}
      \end{array}\right)
\end{equation*}

\begin{equation*}
 F_-(U_R)=\dfrac{A_R}{2\sqrt{3T_R}}
\left(\begin{array}{c}
       \dfrac{{m_+}^2-{m_-}^2}{2}\\
      \dfrac{{m_+}^3-{m_-}^3}{3}
      \end{array}\right)
\end{equation*}
where
\begin{equation*}
\left\{\begin{array}{l}
        M_+=\max(0,u_L+\sqrt{3T_L})\\
	M_-=\max(0,u_L-\sqrt{3T_L})\\
	m_+=\min(0,u_R+\sqrt{3T_R})\\
	m_-=\min(0,u_R-\sqrt{3T_R})
       \end{array}\right. 
\end{equation*}
with
\begin{equation*}
 T_L=\dfrac{k}{3\rho \sqrt{\pi}}\sqrt{A_L} \quad \text{and} \quad T_R=\dfrac{k}{3\rho \sqrt{\pi}}\sqrt{A_R}.
\end{equation*}

\paragraph{CFL condition}\label{cfl}
We have to impose a CFL (Courant, Friedrichs, Levy) condition on the timestep to prevent blow up of the numerical values
 and to ensure the positivity of $A$. This classical stability condition writes
\begin{equation*}
\Delta{t}\leq n_{CFL} \dfrac{\min\limits_{i}(\Delta{x_i})}{\max\limits_{i}(|u_i|+c_i)},
\end{equation*}
where $c_i=\sqrt{k\sqrt{A_i}/(2\rho \sqrt{\pi})}$ and $n_{CFL}=1$ (resp. $n_{CFL}=0.5$) for the first order scheme (resp. for
 the second order scheme see \ref{SECTION2o}).

We have to notice that the kinetic flux needs a particular CFL condition (see
 \cite[p.24]{Audusse99})
\begin{equation*}
\Delta{t}\leq n_{CFL}\dfrac{\min\limits_{i}(\Delta{x_i})}{\max\limits_{i}(|u_i|+\sqrt{3T_i})}.
\end{equation*}

The Rusanov is the most diffusive but the most simple to implement, the kinetic one is slightly less diffusive but more cpu
 time consuming. VFRoe-ncv and HLL are comparable, but the first one is a little more time consuming due to entropy correction.
 All these fluxes are compared in  \cite{Delestre10b} in Saint-Venant framework.

\subsection{Boundary conditions}
\subsubsection{Characteristic system}
Of course boundary conditions are very important for artery flow. We will not too much insist on them, and we will not present
 for example the classical Windkessel model or variations. Nevertheless, adapting \cite{Bristeau01} to blood flow, we can write
 the homogeneous form of \eqref{compact-form} under the following non conservative form
\begin{equation}
\left\{\begin{array}{l}
\partial_t A+u\partial_x A+A\partial_x u=0\\
\partial_t u+\dfrac{k}{2\rho \sqrt{\pi}\sqrt{A}}\partial_x A+u\partial_x u=0
\end{array}\right..\label{non-conservative}
\end{equation}
Thanks to the Moens Korteweg velocity rewritten as \(A=c^4\pi (2\rho/k)^2\), we get
\begin{equation}
 \left\{\begin{array}{l}
	  \partial_t \left[4c-u\right]+(u-c)\partial_x \left[4c-u\right]=0\\
	  \partial_t \left[4c+u\right]+(u+c)\partial_x \left[4c+u\right]=0
        \end{array}\right.,
\end{equation}
thus we have \(\dfrac{d(4c-u)}{dt}=0\) (respectively \(\dfrac{d(4c+u)}{dt}=0\)) or the Riemann invariant \(4c-u=Cst\) (resp.
 \(4c+u=Cst\)) along the characteristic curve \(C_-\) (resp. \(C_+\)) defined by the equation \(\dfrac{dx}{dt}=u-c\)
 (resp. \(\dfrac{dx}{dt}=u+c\)).

The boundary conditions will be prescribed thanks to the method of characteristics. We set \(U_{bound}=U(0)\) or \(U(L)\) and
 \(U_{near}=U(\Delta x)\) or \(U(L-\Delta x)\) depending on the considered boundary (\(x=0\) or \(x=L\)).

\subsubsection{Subcritical flow}
 We write here 
 boundary condition for subcritical flow: at the boundary the flow is subcritical if we have \(|u_{bound}|<c_{bound}\) or
 equivalently
\begin{equation}
 \left(u_{bound}-c_{bound}\right)\left(u_{bound}+c_{bound}\right)<0. \label{subcritical-cond}
\end{equation}
Of course it seems to be always the case in blood flows, where $|Q/(Ac)|$ is less than 10 \% in physicological cases. This
 concept is more relevant to Saint-Venant equation, where supercritical flow are easily obtained. 
Two cases are possible: we are given either the cross section \(A\) (the pressure thanks to \eqref{pressure}) or
 the discharge \(Q\).
\begin{itemize}
  \item Given cross section \(A_{bound}=Cst\):

At \(x=0\) the Riemann invariant is constant along the ingoing characteristic \(u-c\), thus we have
\begin{equation}
 \left\{\begin{array}{l}
         A_{bound}=Cst\\
	u_{bound}=u_{near}-4\left[c_{near}-c_{bound}\right]
        \end{array}\right..
\end{equation}
At \(x=L\) the Riemann invariant is constant along the outgoing characteristic \(u+c\), thus we have
\begin{equation}
 \left\{\begin{array}{l}
         A_{bound}=Cst\\
	u_{bound}=u_{near}+4\left[c_{near}-c_{bound}\right]
        \end{array}\right..
\end{equation}
If \eqref{subcritical-cond} is not verified by the couple (\(A_{bound},u_{bound}\)), the flow is in fact supercritical.
  \item Given discharge \(Q_{bound}=Cst\):

Depending on the boundary considered we have either at \(x=0\)
\begin{equation}
 0=-Q_{bound}+(u_{near}-4c_{near})A_{bound}+4c_{bound}A_{bound},\label{left-bound}
\end{equation}
or at \(x=L\)
\begin{equation}
 0=-Q_{bound}+(u_{near}+4c_{near})A_{bound}-4c_{bound}A_{bound}.\label{right-bound}
\end{equation}
We recall that \(c_{bound}\) depends on \(A_{bound}\), we solve \eqref{left-bound} or \eqref{right-bound} to get \(A_{bound}\).
\end{itemize}
The source terms might be considered negligible in the characteristic method.

\subsubsection{Supercritical flow}
As previously mentioned, this is not really a relevant case for flows in arteries. But, in veins, or in collapsible tubes, it
 may be relevant (see Pedley \cite{Pedley80}).
A supercritical inflow corresponds to both \(A_{bound}=Cst_1\) and \(Q_{bound}=Cst_2\) that have to be imposed. A supercritical
 outflow is such that  we have \(A_{bound}=A_{near}\) and \(Q_{bound}=Q_{near}\).

Again, source terms are neglected in the characteristic method.

In order to impose the discharge \(Q_{bound}=A_{bound}u_{bound}\), we can use an other method: the use of the first component
 of the numerical flux is an accurate way to proceed, even if \({\cal F}_1(U_L,U_R)=Q_{bound}(t)\) has no unique solution (and
 the complexity of the problem depends on the numerical flux).

\subsection{Second order scheme}\label{SECTION2o}
In order to get a second order scheme, we use the following algorithm:

\paragraph{Step 1} In order to get second order reconstructed variables \(\left(U_{\bullet\pm},\sqrt{A_0}_{\bullet\pm}\right)\)
 we perform a linear reconstruction on variables \(u\), \(A\), \(\sqrt{A}-\sqrt{A_0}\) then we get the reconstruction on
 \(\sqrt{A_0}\).
 
We consider the scalar function \(s\in \mathbb{R}\), its linear reconstruction is defined by
\begin{equation*}
 s_{i-1/2+}=s_i-\dfrac{\Delta x}{2}Ds_i \quad \text{and}\quad s_{i+1/2-}=s_i+\dfrac{\Delta x}{2}Ds_i
\end{equation*}
where \(D\) is one of the following reconstruction operator \(D_{muscl}\), \(D_{eno}\) and \(D_{enom}\). To get the
 reconstruction operator $D$, we introduce the minmod slope limiter
\begin{equation*}
 \text{minmod}(x,y)=\left\{\begin{array}{ll}
                      \min(x,y) & \text{if}\; x,y\geq 0\\
		      \max(x,y) & \text{if}\; x,y\leq 0\\
			0 & \text{else}
                     \end{array}\right..
\end{equation*}
Some other slope limiters might be found in \cite[p.111-112]{LeVeque02} (the Monotonized Central difference limiter and the
 Superbee limiter).

\subparagraph{MUSCL} With the operator
\begin{equation*}
 D_{muscl} s_i=\text{minmod}\left(\dfrac{s_i-s_{i-1}}{\Delta x},\dfrac{s_{i+1}-s_i}{\Delta x}\right),
\end{equation*}
we get the MUSCL linear reconstruction (Monotonic Upwind Scheme for Conservation Law \cite{vanLeer79}). We can find a
 less diffusive but more complicated form of the MUSCL reconstruction in \cite{Delestre10}.

\subparagraph{ENO} With the operator
\begin{equation*}
 D_{eno} s_i=\text{minmod}\left(\dfrac{s_i-s_{i-1}}{\Delta x}+\theta_{eno} \dfrac{\Delta x}{2}D^2 s_{i-1/2},
\dfrac{s_{i+1}-s_i}{\Delta x}+\theta_{eno} \dfrac{\Delta x}{2}D^2 s_{i+1/2} \right),
\end{equation*}
where
\begin{equation*}
 D^2 s_{i+1/2}=\text{minmod}\left(\dfrac{s_{i+1}-2s_i+s_{i-1}}{\Delta x^2},
\dfrac{s_{i+2}-2s_{i+1}+s_i}{\Delta x^2}\right)
\end{equation*}
and
\begin{equation*}
 \theta_{eno}=1,
\end{equation*}
we perform the ENO linear reconstruction (Essentially Non Oscillatory \cite{Harten86,Harten87}, \cite{Shu88}). This reconstruction is
 more accurate than the MUSCL reconstruction but less stable. In order to get a more stable reconstruction we can take
 \(\theta_{eno}\) in \([0,1]\) (it is recommended to take \(\theta_{eno}=0.25\) \cite[p.236]{Bouchut07}).

\subparagraph{modified ENO} With the operator
\begin{equation*}
 D_{enom}s_i=\text{minmod}\left(D_{eno}s_i,2\theta_{enom}D_{muscl}s_i\right), \;\text{with}\;\theta_{enom}\in[0,1],
\end{equation*}
we get the modified ENO linear reconstruction \cite[p.55-56]{Bouchut04} and \cite[p.235-237]{Bouchut07}.

The previous linear reconstructions applied to \(A\) writes
\begin{equation*}
 A_{i-1/2+}=A_i-\dfrac{\Delta x}{2}DA_i,\; A_{i+1/2-}=A_i+\dfrac{\Delta x}{2}DA_i,
\end{equation*}
thus we have the following conservation relation
\begin{equation*}
 \dfrac{A_{i-1/2+}+A_{i+1/2-}}{2}=A_i.
\end{equation*}
In order to get the same kind of relation on the discharges, the following modified reconstruction is performed on the velocity
 variables \(u\)
\begin{equation*}
 u_{i-1/2+}=u_i-\dfrac{A_{i+1/2-}}{A_i}\dfrac{\Delta x}{2}Du_i,\;
 u_{i+1/2-}=u_i+\dfrac{A_{i-1/2+}}{A_i}\dfrac{\Delta x}{2}Du_i,
\end{equation*}
this allows to get the conservation relation
\begin{equation*}
 \dfrac{A_{i-1/2+}u_{i-1/2+}+A_{i+1/2-}u_{i+1/2-}}{2}=A_i u_i.
\end{equation*}
We perform on \(\Psi=\sqrt{A}-\sqrt{A_0}\) the same linear reconstruction as on \(A\), we get \(\Psi_{\bullet\pm}\) and at
 last we have
\begin{equation*}
 \sqrt{A_{0}}_{i+1/2\pm}=\sqrt{A_{i+1/2\pm}}-\Psi_{i+1/2\pm}.
\end{equation*}

\paragraph{Step 2} From variables \(U_{\bullet -}\) and \(U_{\bullet +}\) we get the following reconstructed variables 
\begin{equation}
 \left\{\begin{array}{l}
         \sqrt{A_{i+1/2L}}=\max(\Psi_{i+1/2-}+\min(\sqrt{A_{0}}_{i+1/2-},\sqrt{A_{0}}_{i+1/2+}),0) \\
	U_{i+1/2L}=(A_{i+1/2L},A_{i+1/2L}u_{i+1/2-})^t \\
	 \sqrt{A_{i+1/2R}}=\max(\Psi_{i+1/2+}+\min(\sqrt{A_{0}}_{i+1/2-},\sqrt{A_{0}}_{i+1/2+}),0) \\
	U_{i+1/2R}=(A_{i+1/2R},A_{i+1/2R}u_{i+1/2+})^t
        \end{array}\right..\label{hydro-order2}
\end{equation}

\paragraph{Step 3} In order to get a consistent well-balanced scheme, it is necessary to add a centered source term \(Fc_i\).
 Scheme \eqref{scheme} becomes
\begin{equation*}
U_i^{n+1}=U_i^n-\dfrac{\Delta t}{\Delta x}\left(F_{i+1/2L}^n-F_{i-1/2R}^n-Fc_i^n\right),
\end{equation*}
where
\begin{equation*}
 \left.\begin{array}{l}
	  F_{i+1/2L}^n=F_{i+1/2}^n+S_{i+1/2-}\\
	  F_{i-1/2R}^n=F_{i-1/2}^n+S_{i-1/2+}
       \end{array}\right.,
\end{equation*}
with
\begin{equation*}
 S_{i+1/2-}=\left(\begin{array}{c}
                   0\\
		  {\cal P}(A_{i+1/2-}^n)-{\cal P}(A_{i+1/2L}^n)
                  \end{array}\right)
\end{equation*}
\begin{equation*}
 S_{i-1/2+}=\left(\begin{array}{c}
                   0\\
		  {\cal P}(A_{i-1/2+}^n)-{\cal P}(A_{i-1/2R}^n)
                  \end{array}\right)
\end{equation*}
and
\begin{equation*}
 Fc_i=\left(\begin{array}{c}
             0\\
	    \dfrac{k}{\rho \sqrt{\pi}} \dfrac{ (A^n_{i+1/2-})^{3/2}-(A^n_{i-1/2+})^{3/2}}{3} \Delta \sqrt{A_0}_i
            \end{array}\right),
\end{equation*}
with $\Delta \sqrt{A_0}_i=\sqrt{A_0}_{i+1/2-}-\sqrt{A_0}_{i-1/2+}$. The numerical flux $F^n_{i+1/2}$ is calculated as at first
 order \eqref{num-flux}, {\it i.e.} with the variables obtained with the hydrostatic reconstruction \eqref{hydro-order2}.

Thus we get a second order scheme in space which writes
\begin{equation*}
 U^{n+1}=U^n+\Delta t \Phi(U^n)\;\text{with}\;U=(U_i)_{i\in\mathbb Z}
\end{equation*}
and
\begin{equation}
 \Phi(U_i^n)=-\dfrac{1}{\Delta x}\left(F_{i+1/2L}^n-F_{i-1/2R}^n-Fc_i^n\right).\label{phi}
\end{equation}

\paragraph{Step 4} The second order in time is recovered by a second order TVD (Total Variation Diminishing)
 Runge-Kutta (see \cite{Shu88}), namely the Heun method
\begin{equation}
\left.\begin{array}{l}
\widetilde U^{n+1}=U^{n}+\Delta{t}\Phi(U^{n}),\\
\widetilde U^{n+2}=\widetilde U^{n+1}+\Delta{t}\Phi(\widetilde U^{n+1}),\\
U^{n+1}=\dfrac{U^n+\widetilde U^{n+2}}{2},
\end{array}\right.\label{Heun}
\end{equation}
where $\Phi$ is defined by \eqref{phi}. This method is a kind of predictor-corrector method.
 
A more complete modelling of the arterial flow will imply other source terms (diffusion, tension,...). To observe those effects
 one has to use a higher order scheme, as performed for shallow water in \cite{Noelle06} and \cite{Marche08}.

\section{Validation} \label{validation}
\subsection{About the chosen examples}
The chosen examples are more for sake of illustration of the previous scheme rather than for validation on real clinical datas.
 So, the extremities of the arteries will be non reflecting and the numerical values are more indicative rather than clinically
 relevant. Those examples will show that the scheme behaves well and that some validations from linearised or exact solutions may be
 reobtained. More specifically we insist on the cases with large change of initial section like sudden constriction and sudden
 expansion. Even we present here an aneurism, flow in a stenosis could also be treated.

\subsection{The ideal tourniquet}
In the Saint-Venant community, the dam break problem is a classical test case \cite{Stoker57} (in compressible gas dynamics, it
 is referred as the Sod tube, LeVeque \cite{LeVeque92} or Lighthill \cite{Lighthill78} or \cite{Ockendon83}). We have an
 analytical solution of this problem  thanks to the characteristic method. So, we consider the analogous of this  problem in
 blood flow: a tourniquet is applied and we remove it instantaneously. Of course, this is done in supposing that the pulsatile
 effects are neglected,  so it is more a test case rather than a real tourniquet. Any way, with the previous hypothesis, we get
 a Riemann problem. This test allows us to study the ability of the numerical scheme to give the front propagation properly,
 we notice that the non linear terms are tested (but neither the viscous ones nor the change of basic radius).
\begin{figure}[htbp]
     \begin{center}
       \subfloat[]{\includegraphics[angle=-90,width=0.48\textwidth]{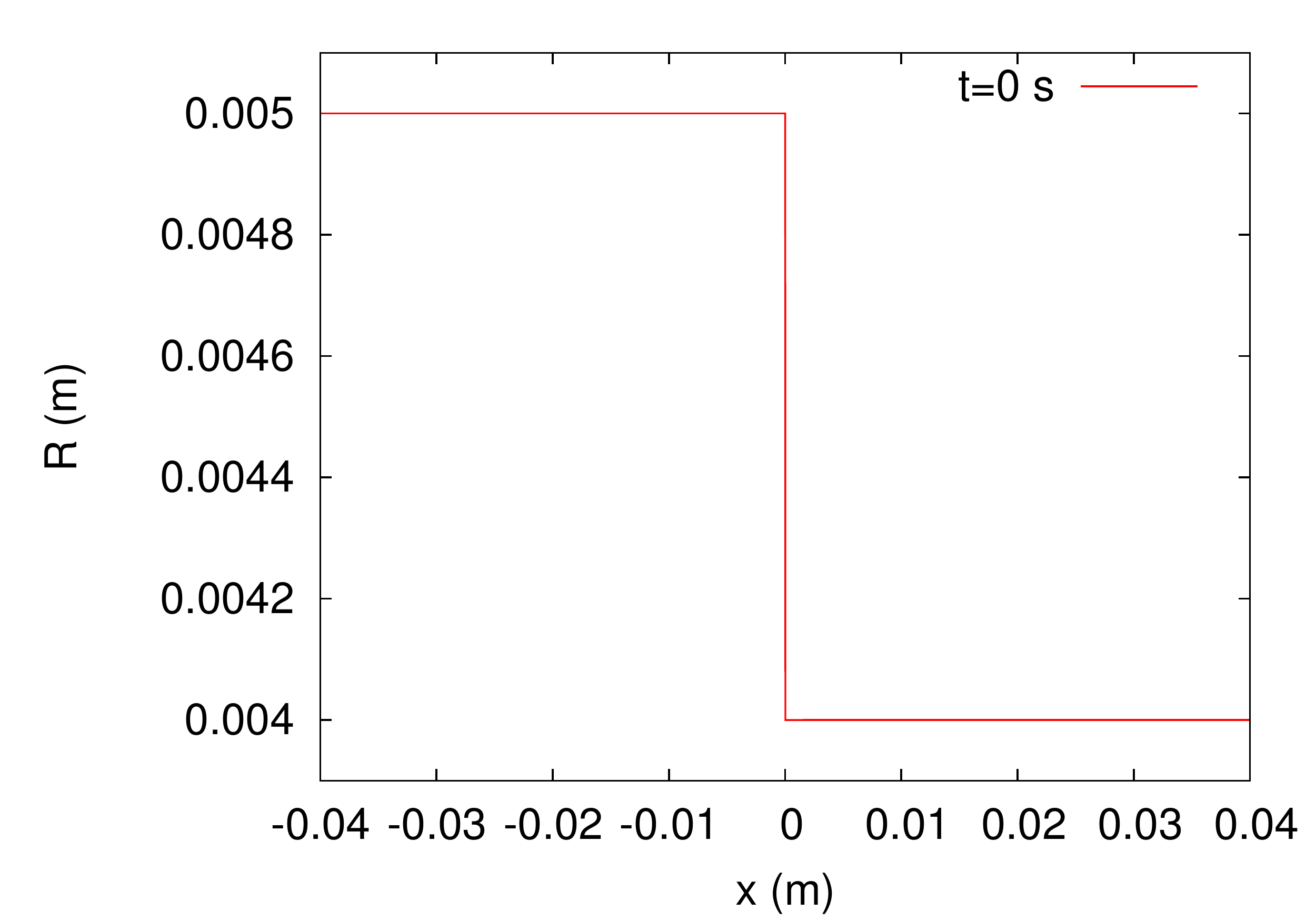}
         \label{fig-tourniquetIniR}}
       \subfloat[]{\includegraphics[angle=-90,width=0.48\textwidth]{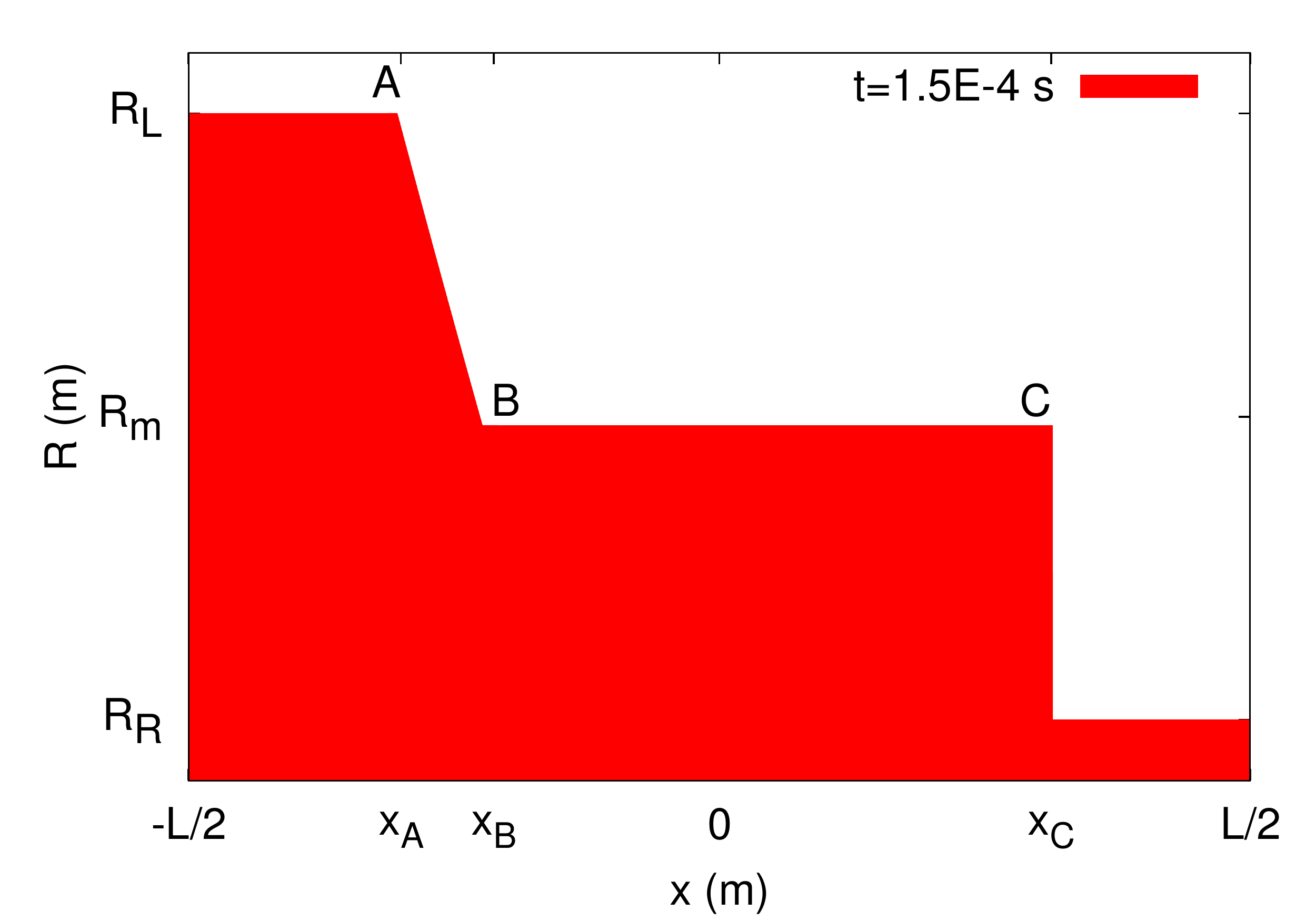}
         \label{fig-tourniquetSolution}}\\
     \end{center}
     \caption{(a) The ideal tourniquet radius initial conditions: a discontinuous initial radius (or pression) is imposed with
 no velocity (the blood is blocked),
 (b) analytical solution of the radius at time $t>0$: a shock wave is moving downstream (to the right) and an expansion wave
 upstream (to the left).}
     \label{fig-tourniquetIni}
\end{figure}
As we consider an artery with a constant radius at rest and without friction, we consider the system \eqref{system2} without
 source term. $L$ is the length of the arteria with $x\in[-L/2,L/2]$. The initial conditions are
\begin{equation}
A(x,0)= \left\{\begin{array}{l}
         A_L\;\text{for}\;x\leq 0\\
         A_R\;\text{for}\;x>0\\
        \end{array}\right.,
\end{equation}
with $A_L>A_R$ and $Q(x,0)=u(x,0)=0$ (see figure \ref{fig-tourniquetIni}).

 As for the dam break, with the characteristic method we get an analytical solution. At time $t>0$, from the upstream to the
 downstream of the artery, we have four different states (as illustrated for the radius $R=\sqrt{A/\pi}$ in Figure
 \ref{fig-tourniquetSolution}):

\begin{enumerate}
 \item Upstream point A located at $x_A(t)=-c_L t$, the state is the same as the initial one:
\begin{equation}
 \left\{\begin{array}{l}
         A(x,t)=A_L\\
	 u(x,t)=0
        \end{array}\right..
\end{equation}
\item Between point A and point B located at $x_B(t)=(u_M-c_M)t=(4c_L-5c_M)t$ (the subscript stands for the intermediate state
 in the following zone), we have
\begin{equation}
 \left\{\begin{array}{l}
         u(x,t)=\dfrac{4}{5}\dfrac{x}{t}+\dfrac{4}{5}c_L \\
	\\
	c(x,t)=-\dfrac{1}{5}\dfrac{x}{t}+\dfrac{4}{5}c_L
        \end{array}\right..
\end{equation}
\item Between point B and C located at $x_C(t)=st$, with the Rankine Hugoniot relation $s=A_M u_M/(A_M-A_R)$.
 We have an intermediate state
\begin{equation}
 \left\{\begin{array}{l}
         A(x,t)=A_M\\
	 u(x,t)=u_M
        \end{array}\right.,
\end{equation}
where the variables $s$, $A_M$, $u_M$ and $Q_M=u_M A_M$ are obtained thanks to the following system
\begin{equation*}
\left\{\begin{array}{l}
        u_M+4c_M=u_L+4c_L\\
	Q_R-Q_M=s(A_R-A_M)\\
	\left(\dfrac{Q_R^2}{A_R}+\dfrac{k}{3\rho\sqrt{\pi}}{A_R}^{3/2}\right)
	-\left(\dfrac{Q_M^2}{A_M}+\dfrac{k}{3\rho\sqrt{\pi}}{A_M}^{3/2}\right)
	=s(Q_R-Q_M)
       \end{array}\right..
\end{equation*}
As $u_R=u_L=Q_R=Q_L=0$, this system reduces to
\begin{equation*}
\left\{\begin{array}{l}
        u_M+4c_M=4c_L\\
	Q_M=s(A_M-A_R)\\
	\left(\dfrac{Q_M^2}{A_M}+\dfrac{k}{3\rho\sqrt{\pi}}{A_M}^{3/2}\right)
	-\dfrac{k}{3\rho\sqrt{\pi}}{A_R}^{3/2}
	=sQ_M
       \end{array}\right..
\end{equation*}

This system is solved iteratively.
\item Downstream point C, the state is the same as the initial one
\begin{equation}
 \left\{\begin{array}{l}
         A(x,t)=A_R\\
	 u(x,t)=0
        \end{array}\right..
\end{equation}
\end{enumerate}

\begin{figure}[htbp]
     \begin{center}
       \subfloat[]{\includegraphics[angle=-90,width=0.48\textwidth]{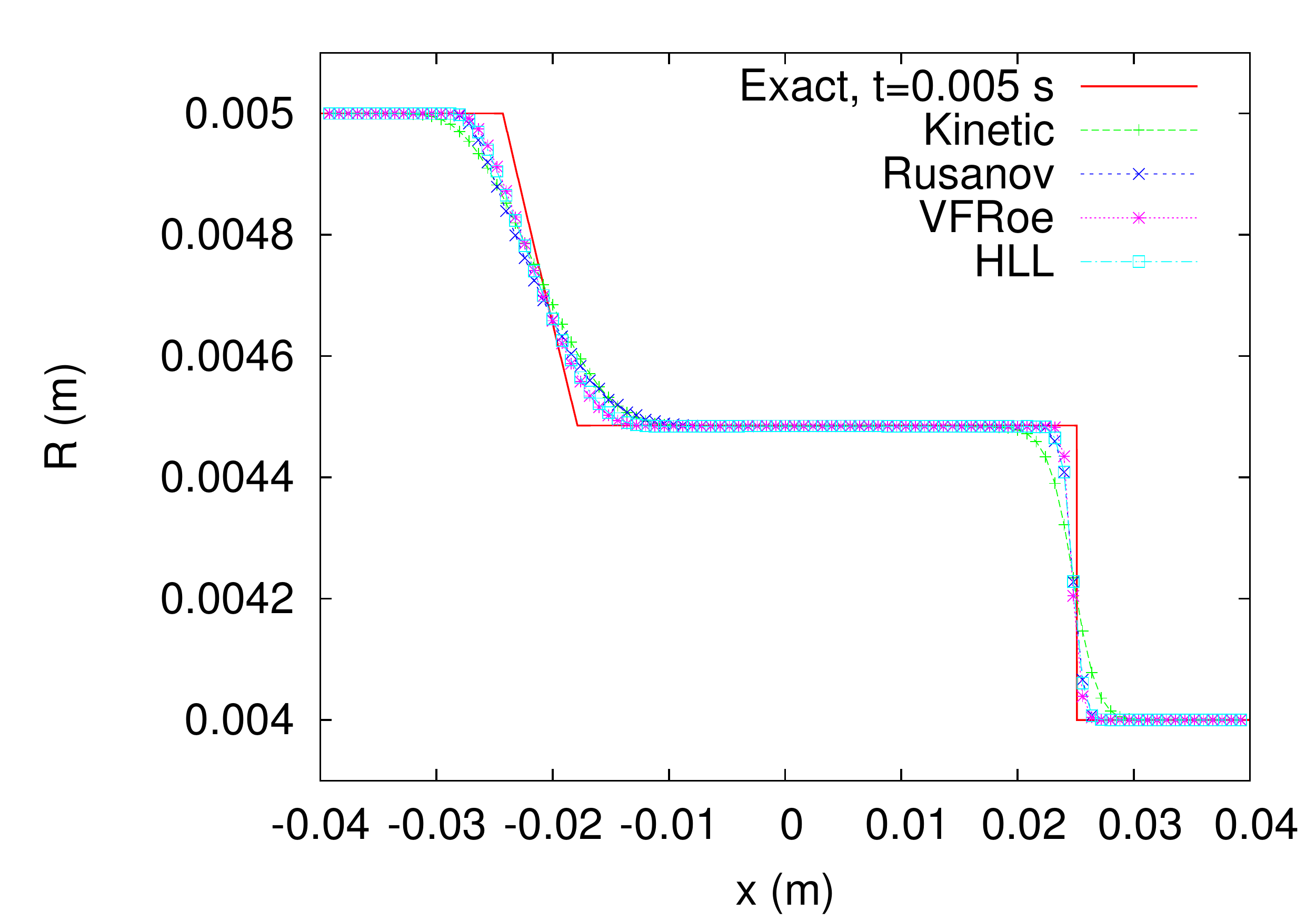}
         \label{fig-tourniquetEndR}}
       \subfloat[]{\includegraphics[angle=-90,width=0.48\textwidth]{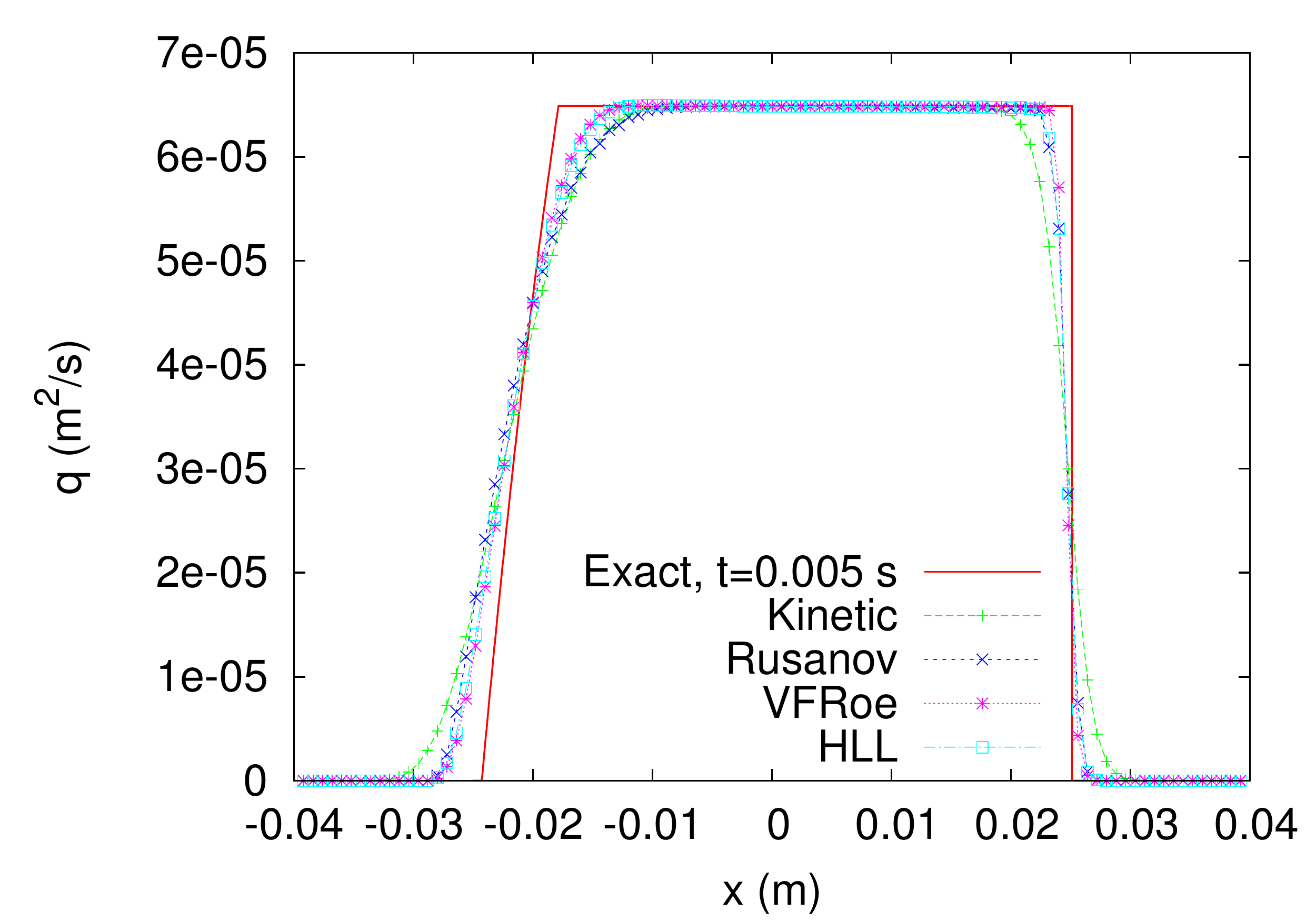}
         \label{fig-tourniquetEndQ}}\\
     \end{center}
     \caption{Snapshots of the ideal tourniquet -- comparison between analytical and numerical results at $t=T_{end}$:
 (a) radius and (b) discharge.}
     \label{fig-tourniquetEnd}
\end{figure}

The simulations were performed with the first order scheme with each of the fluxes, $J=100$ cells and a fix CFL of $1$ ($0.5$
 for the kinetic flux) and we
 have used the following numerical values
 $C_f=0$, $k=1.0\;10^7\text{Pa/m}$, $\rho=1060\text{m}^3$, $R_L=5.\;10^{-3}\text{m}$, $R_R=4.\;10^{-3}\text{m}$,
 $L=0.08 \text{m}$, $T_{end}=0.005 \text{s}$, $c_R=\sqrt{kR_R/(2\rho)}\simeq 4.34\text{m}/\text{s}$,
 $c_L=\sqrt{kR_L/(2\rho)}\simeq 4.86\text{m}/\text{s}$.
An initial flow  at rest: $
 Q(x,t=0)=0 \text{m}^3/\text{s} 
$
and
\begin{equation*}
 A(x,t=0)=\left\{\begin{array}{l}
  \pi {R_L}^2,\;\text{if}\;x\in[-0,04:0] \\
  \pi {R_R}^2,\;\text{if}\;x\in]0:0,04]
 \end{array}\right..
\end{equation*} 
As illustrated in Figure \ref{fig-tourniquetEnd}, we notice that the Rusanov and Kinetic fluxes are more diffusive than the
 others. With mesh refinement and/or scheme order increasing, we should improve the results both at the shock level and the
 expansion wave. As noticed in \cite{Delestre10b} for the shallow water equations, the VFRoe-ncv needs more cpu time than the
 HLL flux. So in what follows we will focus on the HLL flux.
 This is a rather severe test case as in practice,  $|Q/(Ac)|$ is less than 10 \% in physicological cases. 
 
\subsection{Wave equation}
\label{we}

\begin{figure}[htbp]
     \begin{center}
       \subfloat[]{\includegraphics[width=0.48\textwidth]{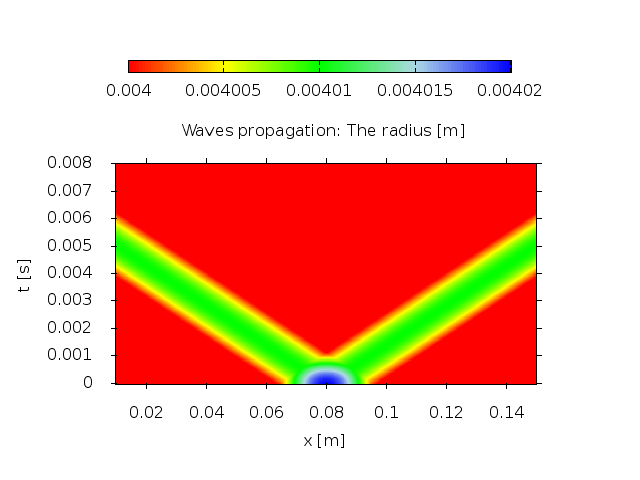}
         \label{fig-Rwave-timeCol}}
       \subfloat[]{\includegraphics[width=0.48\textwidth]{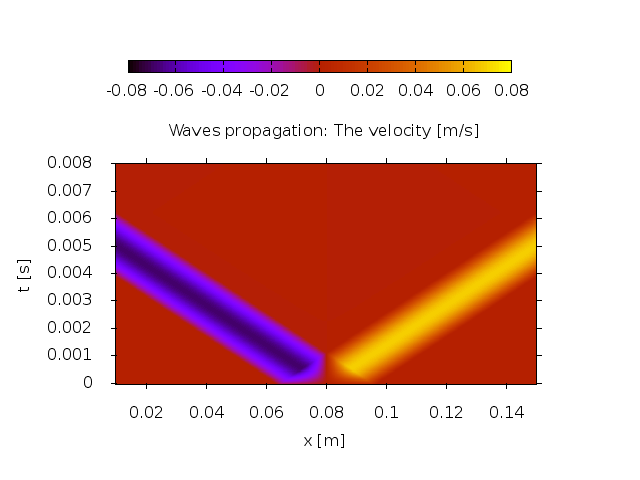}
         \label{fig-Uwave-timeCol}}\\
     \end{center}
     \caption{The propagation of an initial (non physiological) pulse given at $t=0$ in space and time: (a) radius and (b)
 velocity: the initial perturbation splits into a forward and backward travelling wave (colors online).}
     \label{fig-wave-timeCol}
\end{figure}
\begin{figure}[htbp]
     \begin{center}
       \subfloat[]{\includegraphics[angle=-90,width=0.48\textwidth]{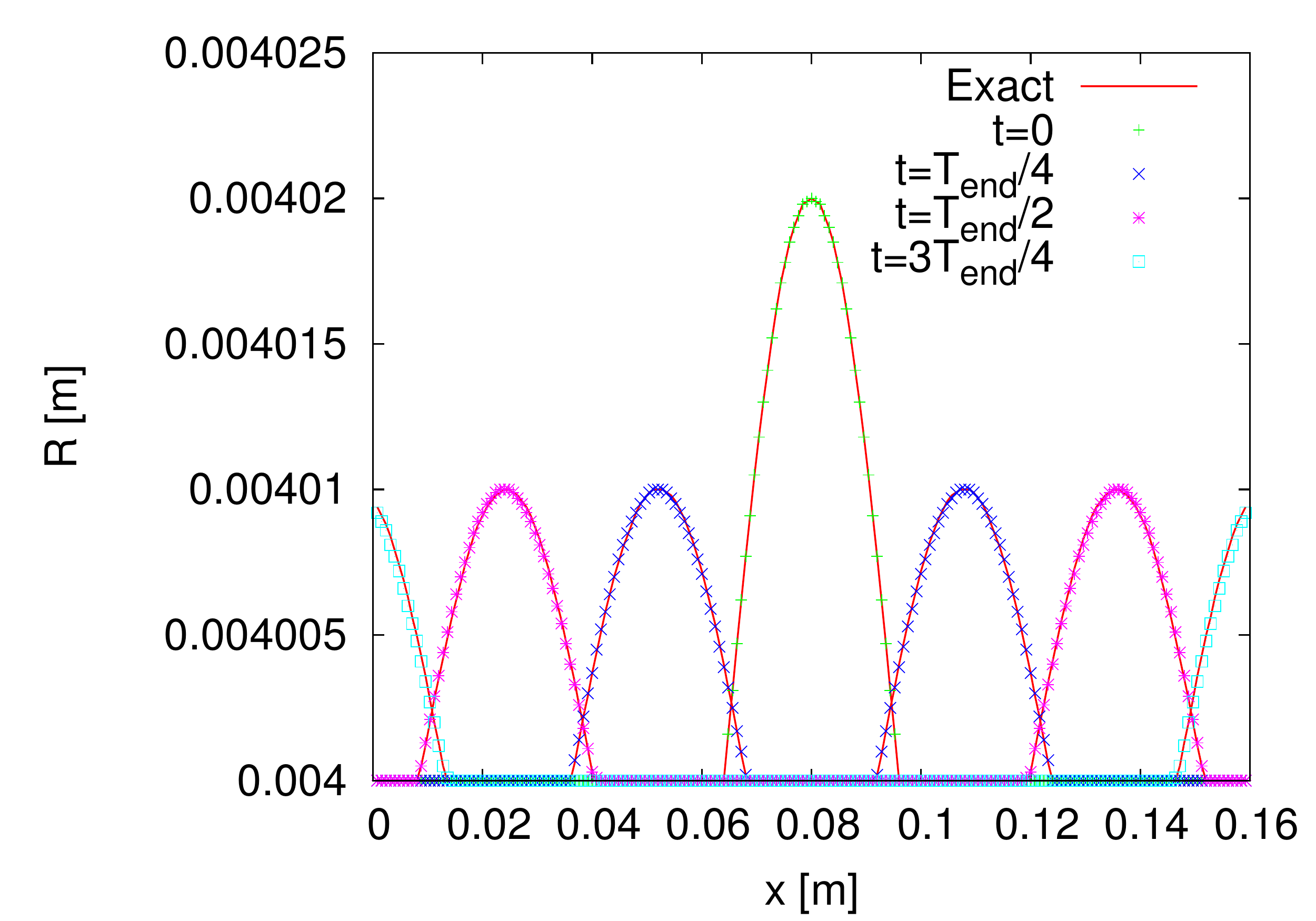}
         \label{fig-waveR}}
       \subfloat[]{\includegraphics[angle=-90,width=0.48\textwidth]{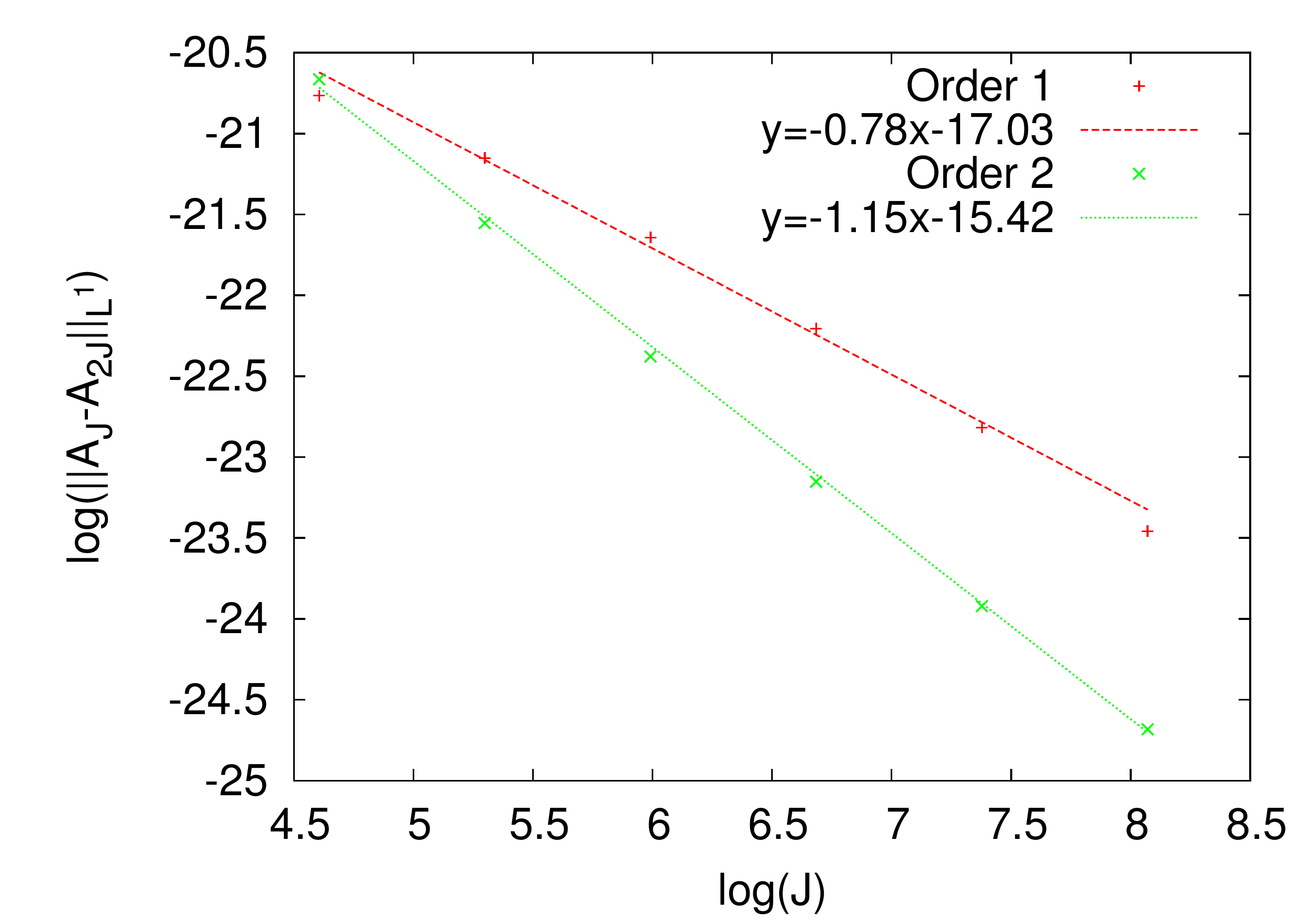}
         \label{fig-waveError}}\\
     \end{center}
     \caption{(a) Radius as function of the position at 3 time steps:
 $t=0$, $t=T_{end}/4$, $t=T_{end}/2$ and $t=3T_{end}/4$, comparison between the first order
 scheme with HLL flux and exact solution and (b) the error made on the area calculation at $t=0.004s$.}
     \label{fig-waveProp}
\end{figure}

So, as in practice, in arteries,  the mean velocity is far smaller than the Moens Korteweg celerity
 $|Q/(Ac)|\ll 1 $, here we will study the ability to catch the relevant small linearised wave propagation speeds. We look at the
 evolution of a given perturbation of radius with no velocity perturbation and see how it evolves. Thus we consider the system
 without friction, with a constant section at rest, a steady state $(R_0,u_0=0)$. Looking at perturbation at a small level
 $\varepsilon \ll1$, we have  $(R,u)=(R_0+\varepsilon R_1,\varepsilon u_1)$ where subscript $1$ refers to the perturbation.
 Substitution in \eqref{system2} and neglecting the small terms, it becomes
the d'Alembert equation for the perturbations $u_1$ and $R_1$:
\begin{equation} 
	  \partial_{t}^2 (u_1,R_1) -{c_0}^2 \partial_{x}^2 (u_1,R_1)=0,
	  \label{lerond}
\end{equation}
the same is valid for variables $Q_1=A_0u_1$ the perturbation of flux and $p_1=kR_1$ the perturbation of pressure
with the Moens Korteweg wave velocity
\begin{equation}
 c_0=\sqrt{\dfrac{k \sqrt{A_0}}{2\rho\sqrt{\pi}}}=\sqrt{\dfrac{k R_0}{2\rho}}.
\end{equation}
With initial conditions $u_1(x,0)=0$, $R_1(x,0)=\phi(x)$ we obtain so analytical solution:
\begin{equation*}
\left\{\begin{array}{l}
 R(x,t)=R_0+\dfrac{\varepsilon}{2}\left[\phi\left(x-c_0 t\right)+\phi\left(x+c_0 t\right)\right] \\
 u(x,t)=-\dfrac{\varepsilon}{2}\dfrac{c_0}{R_0}\left[-\phi\left(x-c_0 t\right)
+\phi\left(x+c_0 t\right)\right]        
       \end{array}\right..
\end{equation*}
The following numerical values have been used for the figure \ref{fig-wave-timeCol}: $J=200$,
 $C_f=0$, $k=1.0\;10^8\text{Pa/m}$, $\rho=1060\text{m}^3$, $R_0=4.\;10^{-3}\text{m}$,
 $L=0.16 \text{m}$, $T_{end}=0.008 \text{s}$, $c_0=\sqrt{kR_0/(2\rho)}\simeq 13.7\text{m}/\text{s}$.
As initial conditions, we take a fluid at rest $Q(x,0)=0 \text{m}^3/\text{s}$ with an initial deformation of the radius
\begin{equation*}
A(x,0)= \left\{\begin{array}{l}
  \pi {R_0}^2,\; \text{if}\;x\in[0:40L/100]\cup [60L/100:L] \\
  \pi {R_0}^2\left[1+\varepsilon \sin\left(\pi \dfrac{x-40L/100}{20L/100}\right)\right]^2 ,\; \text{if}\;x\in]40L/100:60L/100[
 \end{array}\right.,
\end{equation*}
with $\varepsilon=5.10^{-3}$.
As illustrated, in Figure \ref{fig-waveR}, 
we get two waves,  propagating on the  right  (respectively left) with a positive (resp. negative) velocity.
The two waves are represented at several moments. We notice no numerical diffusion.  

\subsection{"The  man at eternal rest"}\label{deadman}
\begin{figure}[htbp]
     \begin{center}
       \subfloat[]{\includegraphics[angle=-90,width=0.48\textwidth]{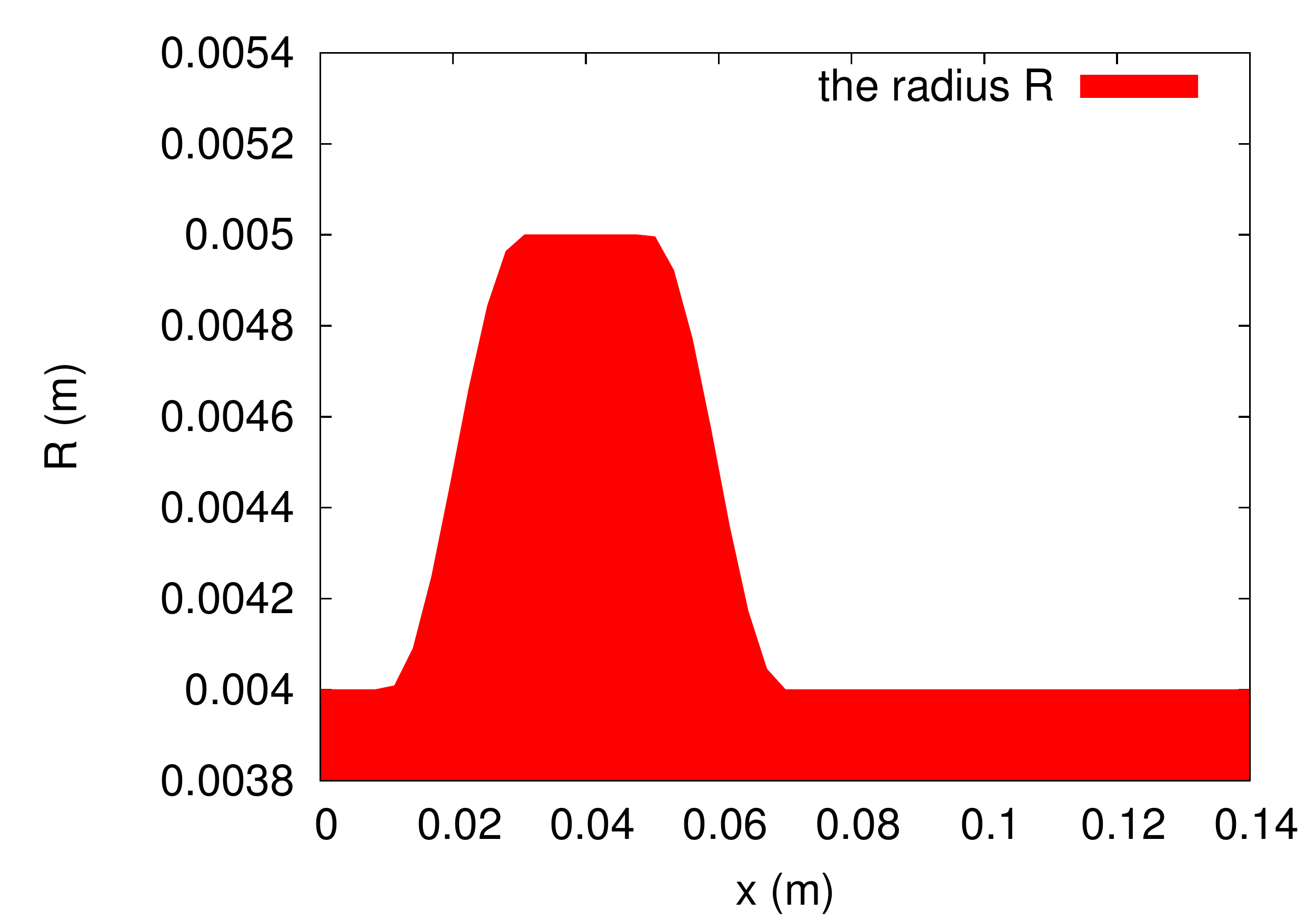}
         \label{fig-deadmanIniR2}}
       \subfloat[]{\includegraphics[angle=-90,width=0.48\textwidth]{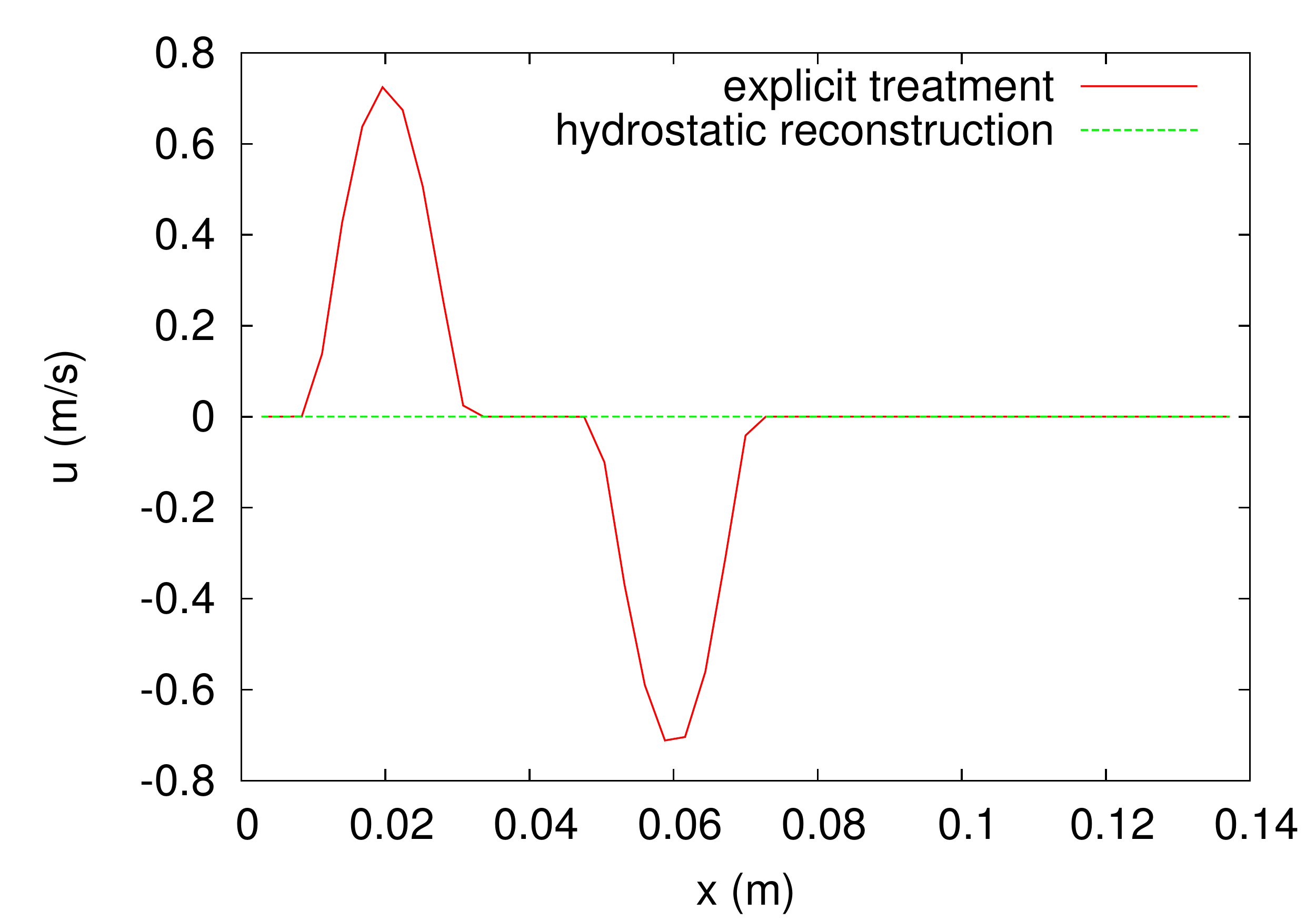}
         \label{fig-source-term-treatment2}}\\
     \end{center}
     \caption{The "dead man case": (a) the radius of the arteria,
 (b) comparison of the velocity at time $t=5\text{s}$ between an explicit treatment of the source term and the hydrostatic
 reconstruction. The spurious flow effect is clearly visible if an inappropriate scheme is chosen.}
     \label{fig-deadman}
\end{figure}

The previous test case did not involve drastic changes in the basic radius. We will show in this section that a not adapted
 source term treatment in \eqref{system2} may give non-physical velocity. 
In this test case, we consider a configuration with no flow and with a change of radius $R_0(x)$, this is the case for a dead
 man with an aneurism. Thus the section of the artery is not constant and the velocity is $u(x,t)=0\;\text{m/s}$.
 We use the following numerical values.
 $J=50$ cells, $C_f=0$, $k=1.0\;10^8\text{Pa/m}$, $\rho=1060\text{m}^3$,
 $L=0.14 \text{m}$, $T_{end}=5 \text{s}$.
As initial conditions, we take a fluid at rest $Q(x,0)=0 \text{m}^3/\text{s}$ and
\begin{equation*}
R(x,0)=R_0(x)= \left\{\begin{array}{l}
  R_o\; \text{if}\;x\in[0:x_1]\cup [x_4:L] \\
  R_o+\dfrac{\Delta R}{2} \left[\sin\left(\left(\dfrac{x-x_1}{x_2-x_1}\pi-\pi/2\right)\right)+1\right] ,\; \text{if}\;x\in]x_1:x_2[\\
R_o+\Delta R\;\text{if}\;x\in[x_2:x_3] \\
R_o+\dfrac{\Delta R}{2}\left(\cos\left(\dfrac{x-x_3}{x_4-x_3}\pi\right)+1\right) ,\; \text{if}\;x\in]x_3:x_4[
 \end{array}\right.,
\end{equation*}
with $R_o=4.0\;10^{-3}\text{m}$, $\Delta R=1.0\;10^{-3}\text{m}$, $x_1=1.0\;10^{-2}\text{m}$, $x_2=3.05\;10^{-2}\text{m}$,
 $x_3=4.95\;10^{-2}\text{m}$ and $x_4=7.0\;10^{-2}\text{m}$ (figure \ref{fig-deadmanIniR2}).

As illustrated in Figure \ref{fig-source-term-treatment2}, with an explicit treatment of the source term in $A_0$, we get non
 zero velocities at the arteria cross section variation level.
 The "man at eternal rest" is not preserved, spurious flow appear (avoiding artifacts such as \cite{Kirkman03}). As expected,
 the hydrostatic reconstruction \eqref{hydro-rec} preserves exactly the steady state at rest.

\subsection{Wave reflection-transmission in an aneurism}
\label{deadwave}
\subsubsection{Analytical linear transmission}
We observed that waves propagate fairly well in a straight tube, we now observe  the reflexion
and transmission trough a sudden change of section (from $A_L$ to $A_R$) in an elastic tube. The d'Alembert equation
 \eqref{lerond} admits harmonic waves $e^{i(\omega t - x/c_0)}$, so that a plane wave is (symbol 1 represents a linear
 perturbation as in section \ref{we}):
$$Q_1=Y_0p_1,\mbox{ with } p_1\propto {\cal{R}}e(e^{i(\omega t - x/c_0)})$$
and $Y_0=A_0/(\rho c_0)$ is called the characteristic admittance. 
With this definition, we can look at a change of section from a radius $R_L$ to a radius $R_R$: 
an imposed right traveling plane wave 
$e^{i(\omega t - x/c_L)}$ moving at celerity $c_L$ in the left media (subscript $L$) will experience a transmission in the
 right media (subscript $R$) (and will move at celerity $c_R$) and a reflexion (and move at celerity $-c_L$). The coefficients
 of transmission and reflexion depend on the two admittances of the two media $L$ and $R$ (see Lighthill or Pedley):
$$Tr=\frac{2 Y_L}{Y_L+Y_R}\mbox{ and } Re=\frac{Y_L-Y_R}{Y_L+Y_R}.$$
As the admittance does not depend on the frenquency $\omega$, any signal (by Fourier decomposition) will be transmitted and
 reflected with those values.

Here we will study the reflection and the transmission of a small wave in an aneurism.

\subsubsection{Propagation of a pulse to and from an expansion}
\begin{figure}[htbp]
     \begin{center}
       \subfloat[]{\includegraphics[width=0.48\textwidth]{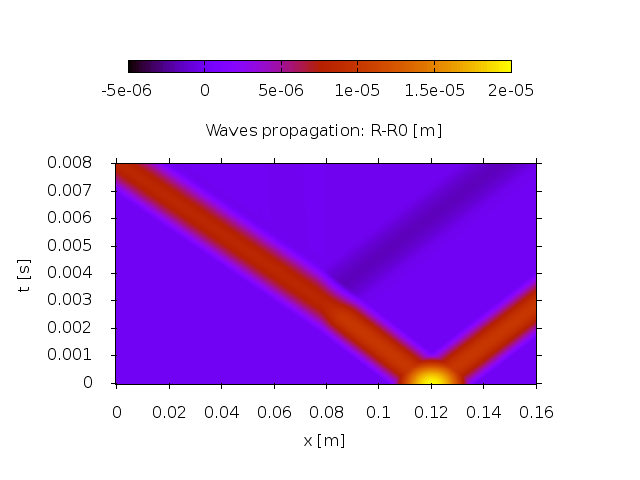}
         \label{fig-Rwave-timeCol1}}
       \subfloat[]{\includegraphics[width=0.48\textwidth]{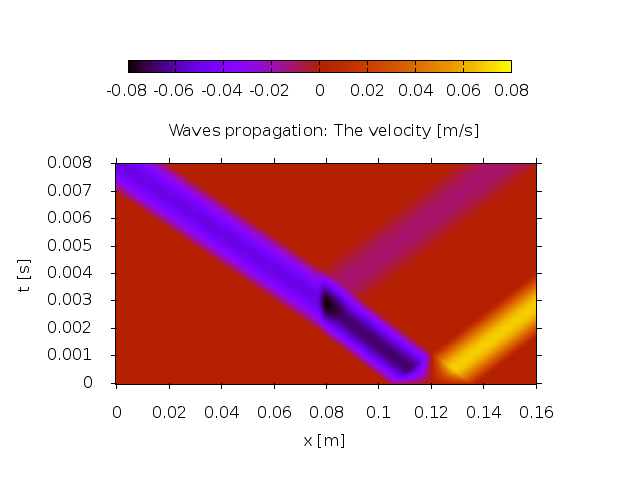}
         \label{fig-Uwave-timeCol1}}\\
     \end{center}
     \caption{To an expansion, the propagation of an initial pulse in space and time: (a) radius and (b) velocity
 (colors online).}
     \label{fig-wave-timeCol1}
\end{figure}

\begin{figure}[htbp]
     \begin{center}
       \includegraphics[angle=-90,width=0.60\textwidth]{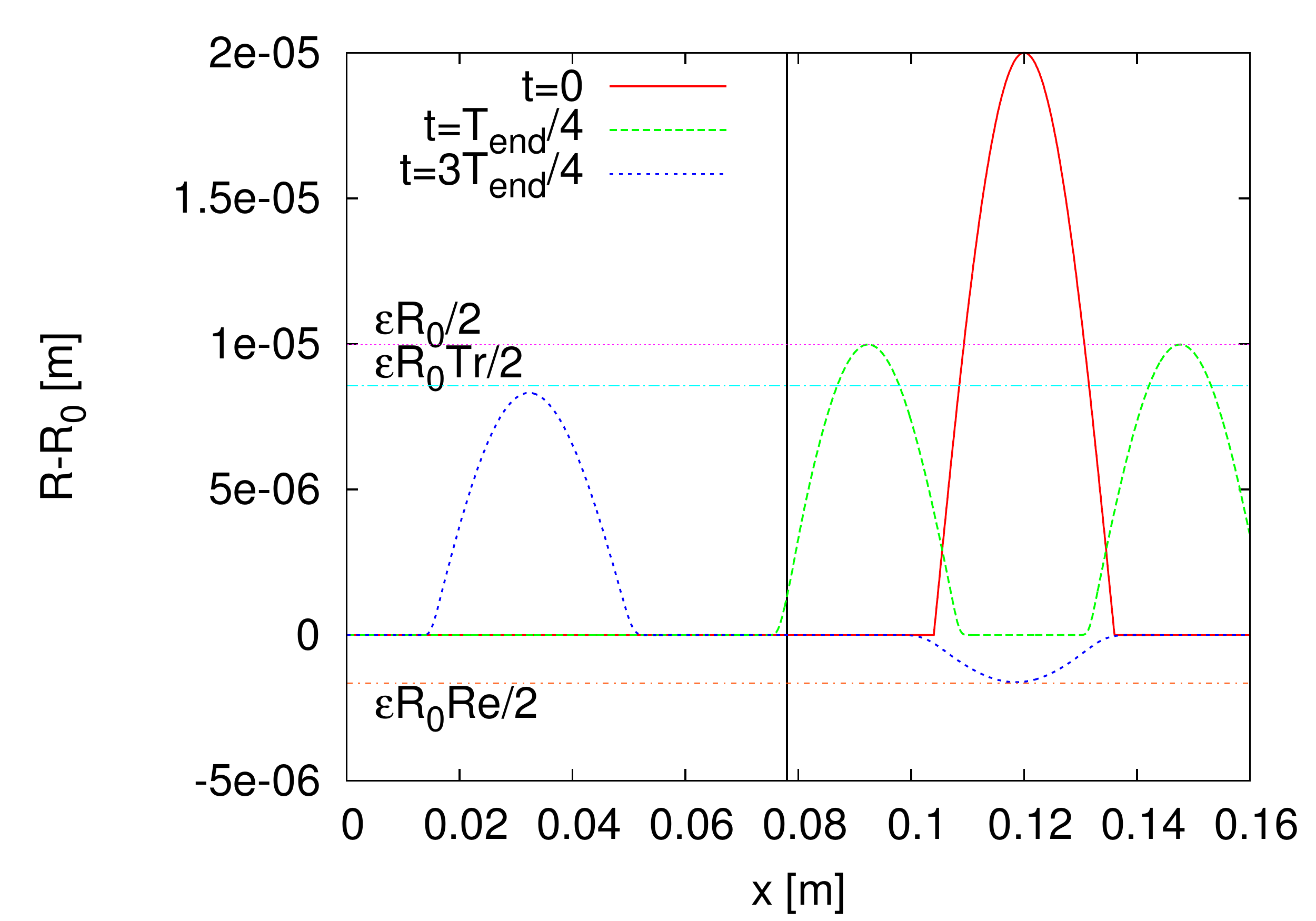}
              \end{center}
     \caption{Radius as function of the position at 3 time steps:
 $t=0$, $t=T_{end}/4$, $t=T_{end}/2$ and $t=3T_{end}/4$.
 The dotted lines represent the level of the predicted reflexion ($Re$) and transmission ($Tr$) coefficients.}
\label{fig-WaveTransReflCol1}
\end{figure}

First we test the case of a pulse in a section $R_R$ passing trough an expansion:  $A_L>A_R$, taking the following numerical
 values:
 $J=1500$ cells, 
 $C_f=0$, $k=1.0\;10^8\text{Pa/m}$, $\rho=1060\text{m}^3$, $R_L=5.\;10^{-3}\text{m}$, $R_R=4.\;10^{-3}\text{m}$,
 $\Delta R=1.0\;10^{-3}\text{m}$, $L=0.16 \text{m}$, $T_{end}=8.0\;10^{-3} \text{s}$, $c_L=\sqrt{kR_L/(2\rho)}\simeq 15.36\text{m}/\text{s}$
 and $c_R=\sqrt{kR_R/(2\rho)}\simeq 13.74\text{m}/\text{s}$.
We take a decreasing shape on a rather small scale:
\begin{equation}
 R_0(x)=\left\{\begin{array}{ll}
                R_R+\Delta R=R_b & \text{if}\;x\in[0:x_1] \\
	      R_R+\dfrac{\Delta R}{2}\left[1+\cos\left(\dfrac{x-x_1}{x_2-x_1}\pi\right)\right] & \text{if}\; x\in]x_1:x_2]\\
	      R_R & \text{else}
               \end{array}\right.,
\end{equation}
with $x_1=19L/40$ and $x_2=L/2$. As initial conditions, we consider a fluid at rest $Q(x,0)=0\text{m}^3/\text{s}$ and the
 following perturbation of radius:
\begin{equation*}
 R(x,0)=\left\{\begin{array}{ll}
                R_0(x)\left[1+\varepsilon \sin\left(\dfrac{100}{20L}\pi\left(x-\dfrac{65L}{100}\right)\right)\right]
 & \text{if} \; x\in[65L/100:85L/100]\\
R_0(x) & \text{else}
               \end{array}\right.,
\end{equation*}
with $\varepsilon=5.0\;10^{-3}$.
 
\begin{figure}[htbp]
     \begin{center}
       \subfloat[]{\includegraphics[width=0.48\textwidth]{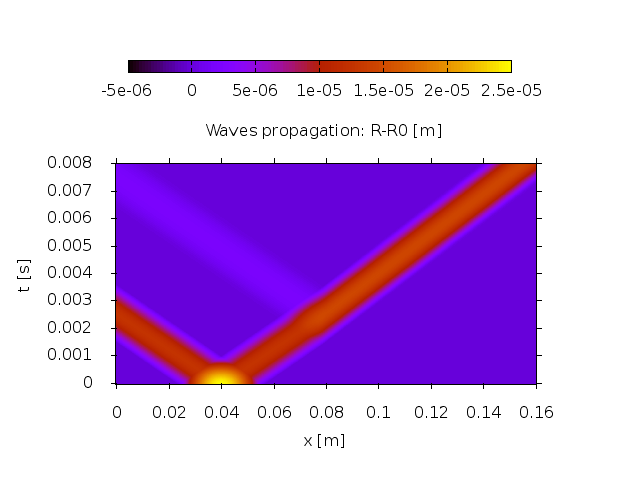}
         \label{fig-Rwave-timeCol2}}
       \subfloat[]{\includegraphics[width=0.48\textwidth]{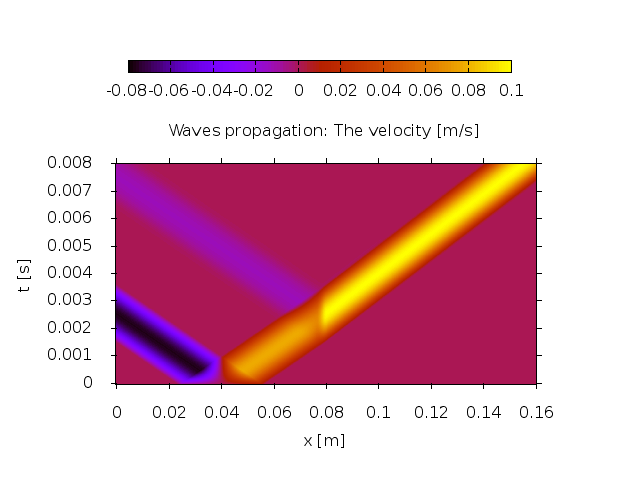}
         \label{fig-Uwave-timeCol2}}\\
     \end{center}
     \caption{From an expansion, the propagation of an initial pulse in space and time: (a) radius and (b) velocity
 (colors online).}
     \label{fig-wave-timeCol2}
\end{figure}

\begin{figure}[htbp]
     \begin{center}
       \includegraphics[angle=-90,width=0.60\textwidth]{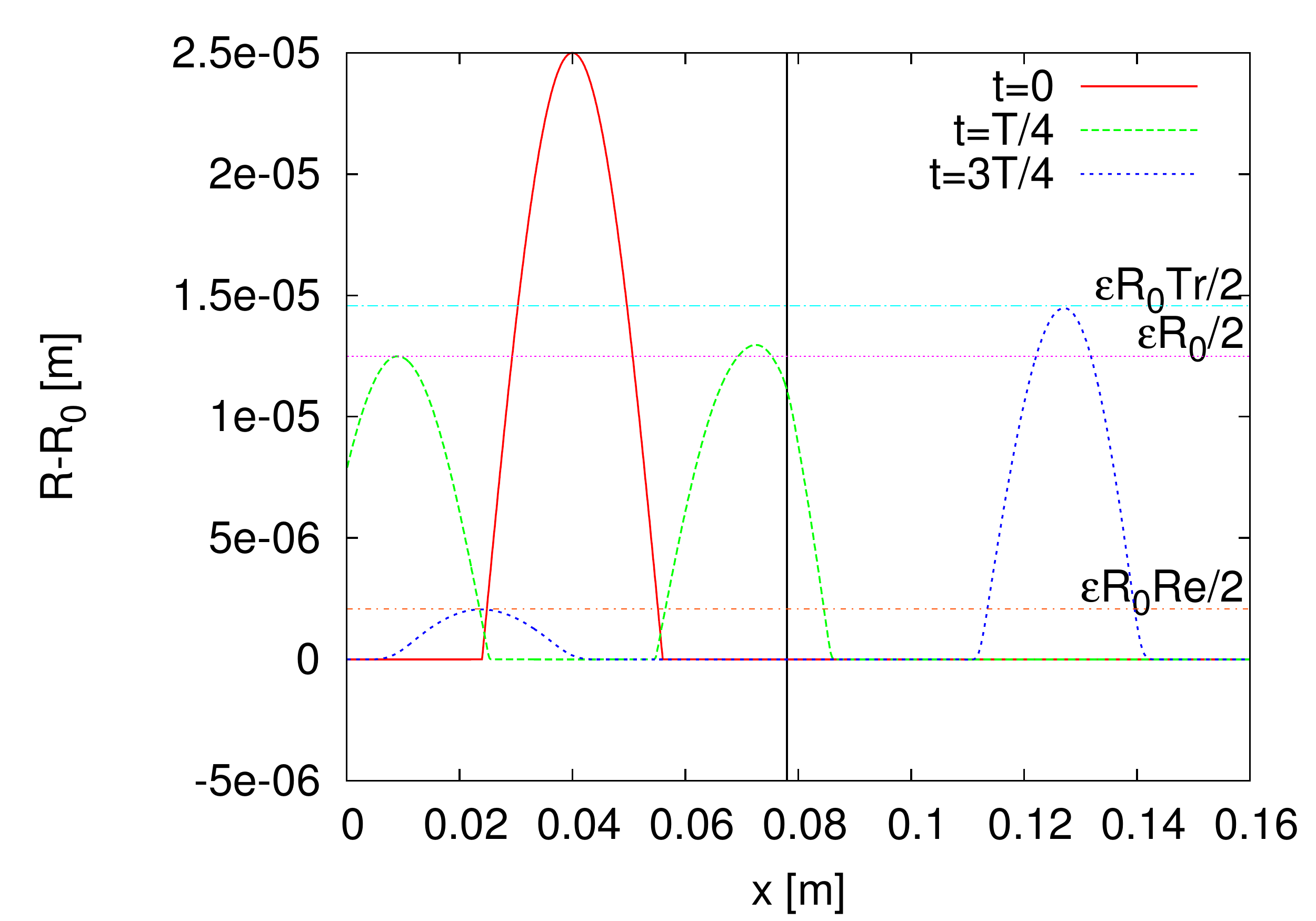}
              \end{center}
     \caption{Radius as function of the position at 3 time steps:
 $t=0$, $t=T_{end}/4$, $t=T_{end}/2$ and $t=3T_{end}/4$.
The dotted lines represent the level of the predicted reflexion ($Re$) and transmission ($Tr$) coefficients. }
\label{fig-WaveTransReflCol2}
\end{figure}

Results showing the wave propagation and expansion are in Figure \ref{fig-wave-timeCol1}. In Figure \ref{fig-WaveTransReflCol1}
 we have the amplitude of the pulse before and after the change of section, the level of the analytical prediction of $T$ and
 $R$ is plotted as well showing that the levels are preserved. 

The same is done for a pulse propagating from an expansion. So, just the radius is changed: 
\begin{equation*}
 R(x,0)=\left\{\begin{array}{ll}
                R_0(x)\left[1+\varepsilon \sin\left(\dfrac{100}{20L}\pi\left(x-\dfrac{15L}{100}\right)\right)\right]
 & \text{if} \; x\in[15L/100:35L/100]\\
R_0(x) & \text{else}
               \end{array}\right.,
\end{equation*}
with $\varepsilon=5.0\;10^{-3}$. Similar results showing the wave propagation are in Figure \ref{fig-Rwave-timeCol2}
and   \ref{fig-Uwave-timeCol2}. In Figure \ref{fig-WaveTransReflCol2} the amplitude of the pulse corresponds to the level
 predicted by the analytical linearised solution.

\subsubsection{The aneurism}

\begin{figure}[htbp]
     \begin{center}
       \subfloat[]{\includegraphics[width=0.48\textwidth]{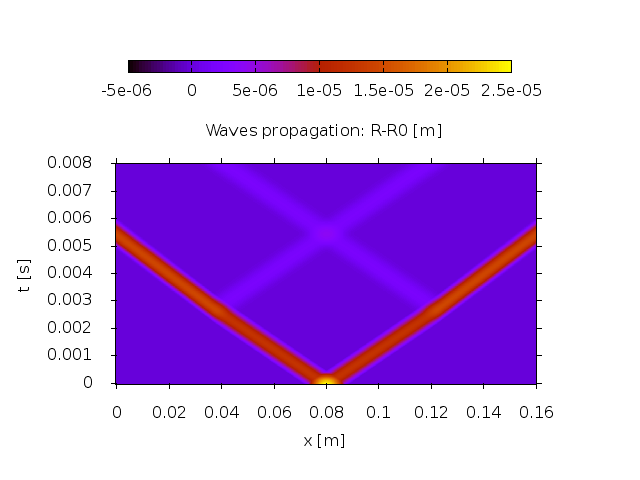}
         \label{fig-Rwave-timeCol3}}
       \subfloat[]{\includegraphics[width=0.48\textwidth]{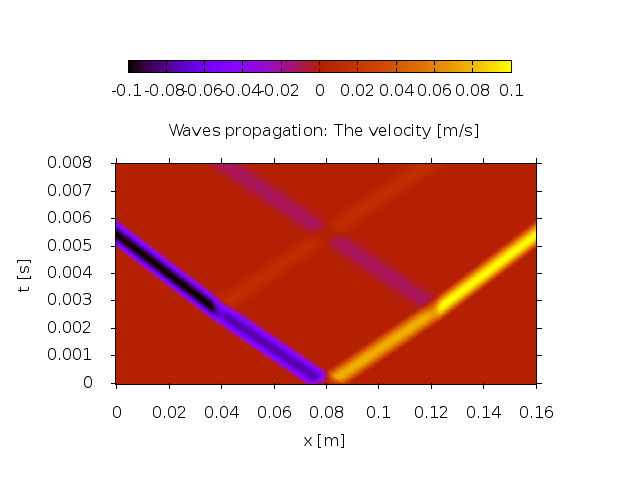}
         \label{fig-Uwave-timeCol3}}\\
     \end{center}
     \caption{The propagation in an aneurysm of an initial pulse in space and time: (a) radius
 and (b) velocity. We note the expected reflexions at the position where the vessel changes his shape (colors online).}
     \label{fig-wave-timeCol3}
\end{figure}

In this part we put an expansion and one constriction modelling an aneurism. An initial pulse at the center $C$ of this aneurism
 propagates in the right and left direction to the boundaries $B$. One sees in Figure \ref{fig-wave-timeCol3} the reflexions
 and transmissions at the position of the variation of radius. These plots have been done with the following numerical values: 
 $J=1500$ cells, 
 $C_f=0$, $k=1.0\;10^8\text{Pa/m}$, $\rho=1060\text{m}^3$, $R_B=4.\;10^{-3}\text{m}$, $R_C=5.\;10^{-3}\text{m}$,
 $\Delta R=1.0\;10^{-3}\text{m}$, $L=0.16 \text{m}$, $T_{end}=8.0\;10^{-3} \text{s}$,
 $c_B=\sqrt{kR_B/(2\rho)}\simeq 13.74\text{m}/\text{s}$
 and $c_C=\sqrt{kR_C/(2\rho)}\simeq 15.36\text{m}/\text{s}$.
With the given shape
\begin{equation}
 R_0(x)=\left\{\begin{array}{ll}
              R_B & \text{if}\;x\in[0:x_1]\bigcup[x_4:L] \\
	      R_B+\dfrac{\Delta R}{2}\left[1-\cos\left(\dfrac{x-x_1}{x_2-x_1}\pi\right)\right] & \text{if}\; x\in]x_1:x_2]\\
	      R_B+\Delta R=R_C & \text{if}\; x\in]x_2:x_3]\\
	      R_B+\dfrac{\Delta R}{2}\left[1+\cos\left(\dfrac{x-x_3}{x_4-x_3}\pi\right)\right] & \text{if}\; x\in]x_3:x_4[
               \end{array}\right.,
\end{equation}
with $x_1=9L/40$, $x_2=L/4$, $x_3=3L/4$ and $x_4=31L/40$. As initial conditions, we consider a fluid at rest
 $Q(x,0)=0\text{m}^3/\text{s}$ and the initial perturbation:
\begin{equation*}
 R(x,0)=\left\{\begin{array}{ll}
                R_0(x)\left[1+\varepsilon \sin\left(\dfrac{100}{10L}\pi\left(x-\dfrac{45L}{100}\right)\right)\right]
 & \text{if} \; x\in[45L/100:55L/100]\\
R_0(x) & \text{else}
               \end{array}\right.,
\end{equation*}
with $\varepsilon=5.0\;10^{-3}$.

\subsubsection{Non adapted treatment of the source terms}
\begin{figure}[]
     \begin{center}
       \subfloat[]{\includegraphics[angle=0,width=0.48\textwidth]{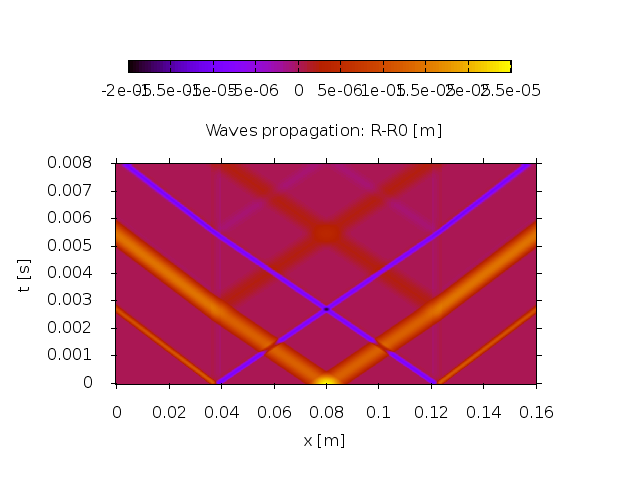}
         \label{fig-deadmanIniR}}
       \subfloat[]{\includegraphics[angle=0,width=0.48\textwidth]{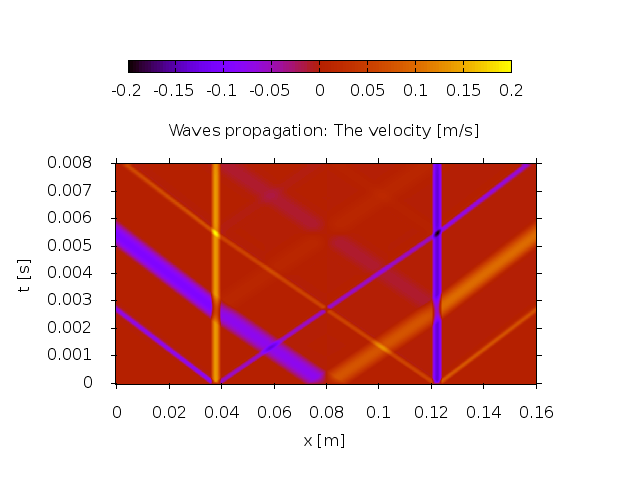}
         \label{fig-source-term-treatment}}\\
     \end{center}
     \caption{A scheme which does not preserve the "man at eternal rest" generates spurious waves were the radius changes
 (colors online).}
     \label{fig-deadman2}
\end{figure}

To insist on the adapted treatment of the source terms, we present in this sub section what happens if a naive scheme is taken
 like in subsection \ref{naive}. We clearly see in Figure \ref{fig-deadman2} that at the change of radius reflected and
 spurious wave appear. The previous constant velocities are present like in Figure \ref{fig-source-term-treatment}. 
Furthermore, the initial data give traveling waves coming from the position of the change of section.
This clearly shows the influence of the source terms in the discretization. 

\subsection{Wave "damping"}\label{wd}

In this last test case, we look at the viscous damping term in the linearized momentum equation. This is the analogous of
 the Womersley \cite{Womersley55} problem, we consider a periodic signal at the inflow which is a perturbation of a steady
 state $(R_0=Cst,u_0=0)$ with a constant section at rest.
We consider the set of equations \ref{system2} with the friction term
under the nonconservative form with the variables $(R,u)$, this system writes
\begin{equation*}
\left\{\begin{array}{l}
	\partial_t R+u \partial_x R+\dfrac{R}{2} \partial_x u=0 \\
        \partial_t u + u\partial_x u + \dfrac{k}{\rho} \partial_x R = - \dfrac{C_f}{\pi} \dfrac{u}{R^2}
       \end{array}\right..
\end{equation*}
We take $R=R_0+\varepsilon R_1$ and $u=0+\varepsilon u_1$, where $(R_1,u_1)$ is the perturbation of the steady state. 
We have after linearisation of the equations:
\begin{equation*}
 \left\{\begin{array}{l}
         2 \partial_t R_1 + R_0 \partial_x u_1=0 \\
	\partial_t u_1 +\dfrac{k}{\rho}\partial_x R_1 = -\dfrac{C_f}{\pi} \dfrac{u_1}{{R_0}^2}
       \end{array}\right.,
\end{equation*}
which may be combined in two wave equations with source term
\begin{equation} 
         \partial_{t}^2 (u_1,R_1) - {c_0}^2 \partial_{x}^2 (u_1,R_1)= -\dfrac{C_f}{\pi}\dfrac{\partial_t (u_1,R_1)}{{R_0}^2}.
        \label{eq:wave-damp}
\end{equation}
Looking for progressive waves ({\it i.e.} under the form $e^{i(\omega t-Kx)}$), we obtain the dispersion relation
\begin{equation}
 K^2=\dfrac{\omega^2}{{c_0}^2}-i\dfrac{\omega C_f}{\pi R_0 {c_0}^2},
 \label{disp}
\end{equation}
with $\omega=2\pi/T_{pulse}$, $T_{pulse}$ the time length of a pulse and $K=k_r+i k_i$, the wave vector not to be confused with
 the stiffness, 
\begin{equation}
 \left\{\begin{array}{l}
         k_r=\left[\dfrac{\omega^4}{{c_0}^4}+ \left(\dfrac{\omega C_f}{\pi {R_0}^2 {c_0}^2}\right)^2\right]^{1/4}
\cos\left(\dfrac{1}{2}\arctan\left(-\dfrac{C_f}{\pi {R_0}^2 \omega}\right)\right) \\
	 k_i=\left[\dfrac{\omega^4}{{c_0}^4}+ \left(\dfrac{\omega C_f}{\pi {R_0}^2 {c_0}^2}\right)^2\right]^{1/4}
\sin\left(\dfrac{1}{2}\arctan\left(-\dfrac{C_f}{\pi {R_0}^2 \omega}\right)\right)
        \end{array}\right..
\end{equation}
\begin{sidewaystable}[htbp]

\begin{center}
{\small
\begin{tabular}{|c|c|c|c|c|c|c|c|c|c|c|}  \cline{1-11}
 \multicolumn{1}{|c|}{$C_f=0.000022$} &
\multicolumn{2}{|c|}{HLL1-SI} &
\multicolumn{2}{|c|}{HLL1-AT} &
\multicolumn{2}{|c|}{HLL MUSCL-SI} &
\multicolumn{2}{|c|}{HLL ENO-SI} &
\multicolumn{2}{|c|}{HLL MUSCL-AT}
 \\ \hline
\multicolumn{1}{|c|}{J} &
\(||Q-Q_{ex}||_{L^1}\) & \(t_{cpu}\) [s]&
\(||Q-Q_{ex}||_{L^1}\) & \(t_{cpu}\) [s]&
\(||Q-Q_{ex}||_{L^1}\) & \(t_{cpu}\) [s]&
\(||Q-Q_{ex}||_{L^1}\) & \(t_{cpu}\) [s]&
\(||Q-Q_{ex}||_{L^1}\) & \(t_{cpu}\) [s]\\ \hline
50  & 0.311E-7 & 0.49   & 0.309E-7 & 0.51   & 0.753E-8 & 2.7    & 0.796E-8 & 2.52   & 0.753E-8 & 2.9  \\ \hline
100 & 0.158E-7 & 1.95   & 0.157E-7 & 2      & 0.232E-8 & 10.65  & 0.216E-8 & 9.95   & 0.232E-8 & 11.47 \\ \hline
200 & 0.789E-8 & 7.8    & 0.783E-8 & 7.96   & 0.121E-8 & 41.95  & 0.129E-8 & 39.3   & 0.121E-8 & 45.31 \\ \hline
400 & 0.386E-8 & 31.26  & 0.383E-8 & 31.75  & 0.458E-9 & 167.16 & 0.651E-9 & 156.47 & 0.457E-9 & 179.4 \\ \hline
800 & 0.183E-8 & 126.06 & 0.182E-8 & 126.92 & 0.267E-9 & 665.91 & 0.401E-9 & 625.62 & 0.267E-9 & 717.64 \\ \hline
Regression & y=-1.02x-13.27 & & y=-1.02x-13.27 & & y=-1.2x-14.19 & & y=-1.04x-14.86 & & y=-1.2x-14.19 &  \\ \hline
 \multicolumn{11}{c}{} \\ \hline

 \multicolumn{1}{|c|}{$C_f=0.000202$} &
\multicolumn{2}{|c|}{HLL1-SI} &
\multicolumn{2}{|c|}{HLL1-AT} &
\multicolumn{2}{|c|}{HLL MUSCL-SI} &
\multicolumn{2}{|c|}{HLL ENO-SI} &
\multicolumn{2}{|c|}{HLL MUSCL-AT}
 \\ \hline
\multicolumn{1}{|c|}{J} &
\(||Q-Q_{ex}||_{L^1}\) & \(t_{cpu}\) [s]&
\(||Q-Q_{ex}||_{L^1}\) & \(t_{cpu}\) [s]&
\(||Q-Q_{ex}||_{L^1}\) & \(t_{cpu}\) [s]&
\(||Q-Q_{ex}||_{L^1}\) & \(t_{cpu}\) [s]&
\(||Q-Q_{ex}||_{L^1}\) & \(t_{cpu}\) [s]\\ \hline
50  & 0.269E-7 & 0.49   & 0.311E-7 & 0.51    & 0.426E-8 & 2.63   & 0.557E-8 & 2.74  & 0.513E-8 & 2.81 \\ \hline
100 & 0.138E-7 & 1.97   & 0.157E-7 & 2       & 0.181E-8 & 10.5   & 0.197E-8 & 10.89 & 0.18E-8 & 11.23 \\ \hline
200 & 0.707E-8 & 7.85   & 0.795E-8 & 7.97    & 0.984E-9 & 41.89  & 0.899E-9 & 43.46 & 0.748E-9 & 44.77 \\ \hline
400 & 0.365E-8 & 31.23  & 0.407E-8 & 31.84   & 0.505E-9 & 167.18 & 0.499E-9 & 173.43& 0.425E-9 & 180.52 \\ \hline
800 & 0.192E-8 & 124.92 & 0.212E-8 & 127.31  & 0.283E-9 & 670.25 & 0.327E-9 & 694.69& 0.296E-9 & 718.3 \\ \hline
Regression & y=-0.95x-13.7 & & y=-0.97x-13.5 & & y=-1.02x-15.24& & y=-1.02x-15.2 &  & y=-1.03x-15.29 &  \\ \hline
 \multicolumn{11}{c}{} \\ \hline

 \multicolumn{1}{|c|}{$C_f=0.005053$} &
\multicolumn{2}{|c|}{HLL1-SI} &
\multicolumn{2}{|c|}{HLL1-AT} &
\multicolumn{2}{|c|}{HLL MUSCL-SI} &
\multicolumn{2}{|c|}{HLL ENO-SI} &
\multicolumn{2}{|c|}{HLL MUSCL-AT}
 \\ \hline
\multicolumn{1}{|c|}{J} &
\(||Q-Q_{ex}||_{L^1}\) & \(t_{cpu}\) [s]&
\(||Q-Q_{ex}||_{L^1}\) & \(t_{cpu}\) [s]&
\(||Q-Q_{ex}||_{L^1}\) & \(t_{cpu}\) [s]&
\(||Q-Q_{ex}||_{L^1}\) & \(t_{cpu}\) [s]&
\(||Q-Q_{ex}||_{L^1}\) & \(t_{cpu}\) [s]\\ \hline
50  & 0.26E-7  & 0.51   & 0.396E-7 & 0.52   & 0.162E-7 & 2.68  & 0.159E-7 & 2.79  & 0.334E-7& 2.81  \\ \hline
100 & 0.146E-7 & 1.98   & 0.215E-7 & 2.03   & 0.916E-8 & 10.56 & 0.899E-8 & 11.03 & 0.176E-7& 11.2 \\ \hline
200 & 0.79E-8  & 7.88   & 0.112E-7 & 8.01   & 0.494E-8 & 41.94 & 0.486E-8 & 43.71 & 0.907E-8& 44.69 \\ \hline
400 & 0.411E-8 & 31.32  & 0.576E-8 & 31.9   & 0.257E-8 & 167.32& 0.253E-8 & 174.25& 0.46E-8 & 178.64 \\ \hline
800 & 0.21E-8  & 125.15 & 0.291E-8 & 127.18 & 0.131E-8 & 669.56& 0.129E-8 & 693.64& 0.231E-8& 714.29 \\ \hline
Regression & y=-0.91x-13.9 & & y=-0.94x-13.33 & & y=-0.91x-14.34 & & y=-0.91x-14.36 & & y=-0.96x-13.43 &  \\ \hline
\end{tabular}}
\end{center}

~\caption{The wave "damping": $L^1$ errors on the discharge $Q$ and CPU times $t_{cpu}$ for
 $C_f=0.000022$, $0.000202$ and $0.005053$. (AT) apparent topgraphy, (SI) semi-implicit and (HLL1) first order scheme
 with HLL flux.}
\label{tbl:TabErDamping2}
\end{sidewaystable}

For the numerical applications, we impose the incoming discharge 
\begin{equation*}
Q_b(t)=Q_{amp}\sin(\omega t)\text{m}^3/\text{s},
\end{equation*}
with $Q_{amp}$ the amplitude of the inflow discharge. We get a damping wave in the domain
\begin{equation}
Q(t,x)=
\left\{\begin{array}{ll}
        0 &,\;\text{if}\;k_r x> \omega t\\
       Q_{amp}\sin\left(\omega t-k_r x\right)e^{k_i x} &,\;\text{if}\;k_r x \leq \omega t
       \end{array}\right., \label{dampingDischarge}
\end{equation}
where $Q_b(t)$ is the discharge imposed at $x=0$.
\begin{figure}[htbp]
     \begin{center}
       \subfloat[]{\includegraphics[angle=-90,width=0.48\textwidth]{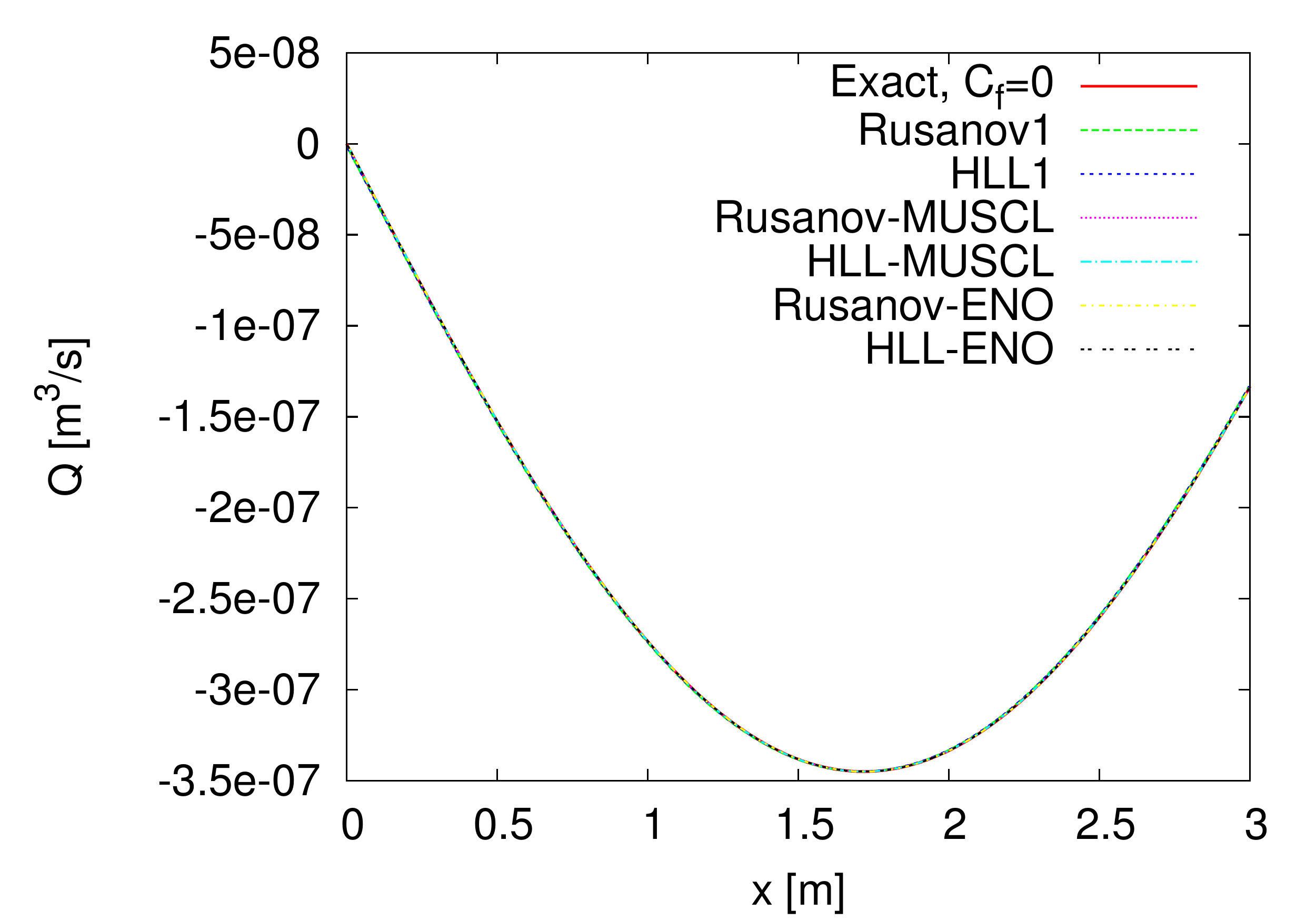}
         \label{fig-wave-damping-Cf0Col}}
       \subfloat[]{\includegraphics[angle=-90,width=0.48\textwidth]{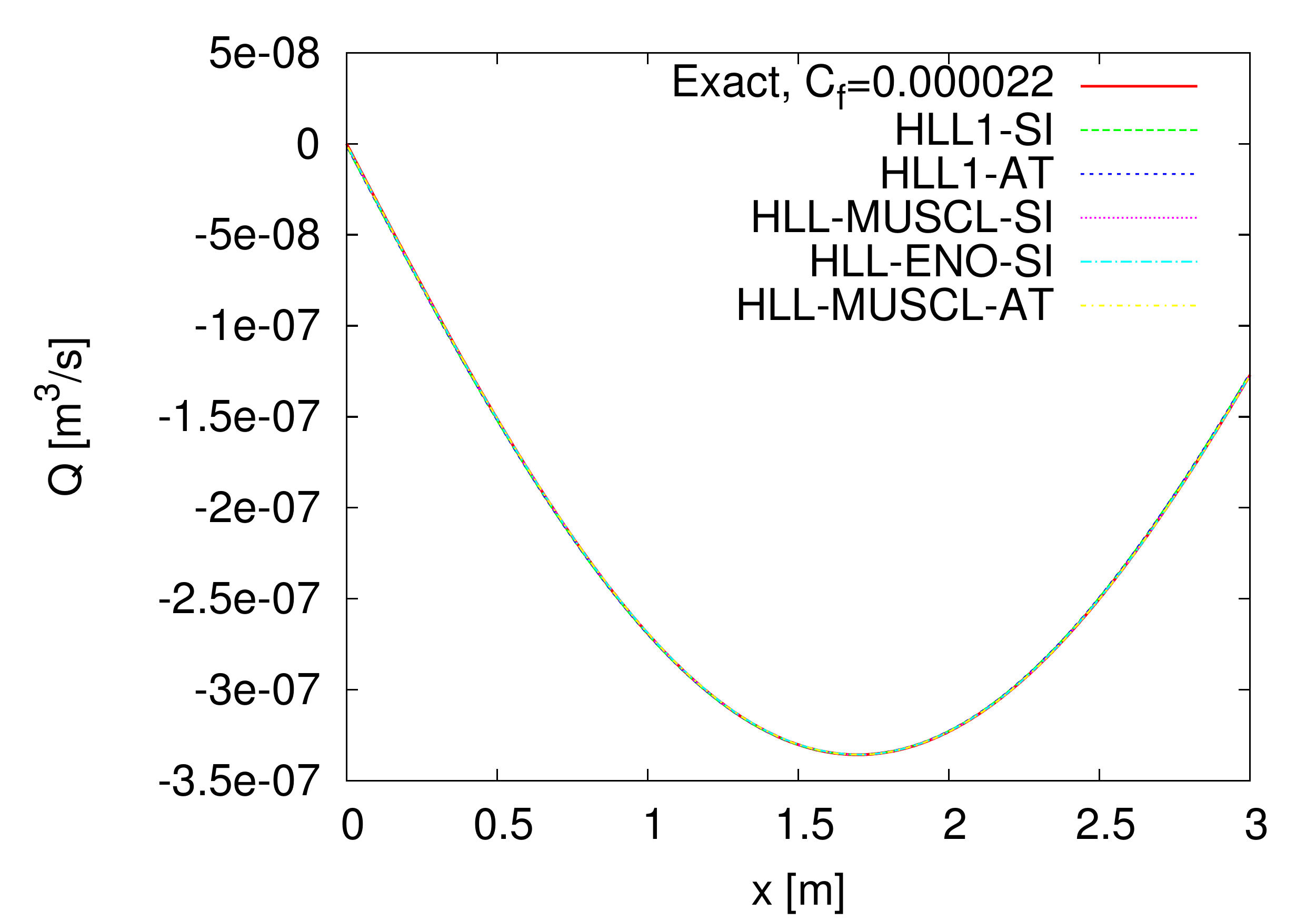}
         \label{fig-wave-damping-Cf0008Col}}\\
       \subfloat[]{\includegraphics[angle=-90,width=0.48\textwidth]{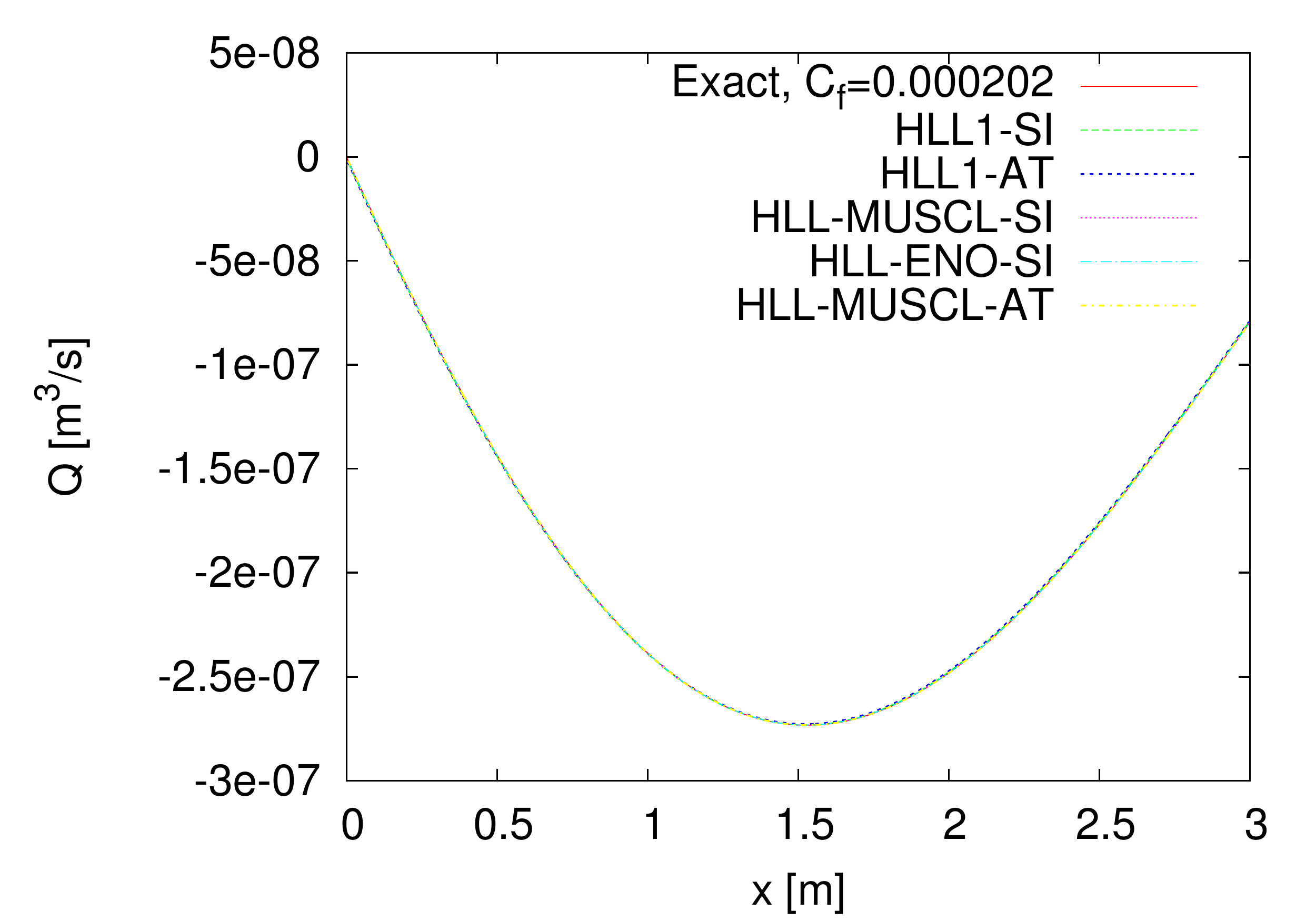}
         \label{fig-wave-damping-Cf0016Col}}
       \subfloat[]{\includegraphics[angle=-90,width=0.48\textwidth]{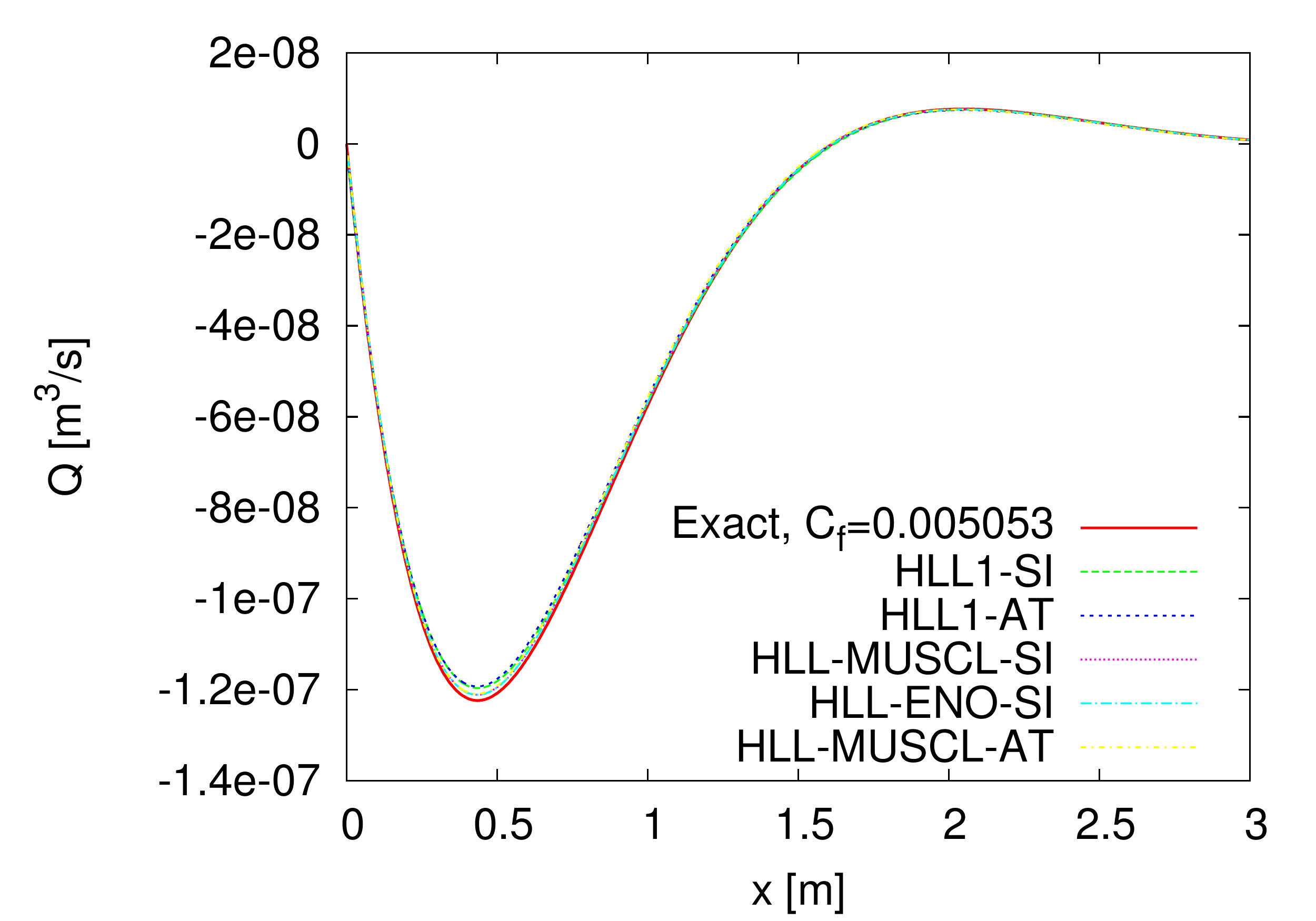}
         \label{fig-wave-damping-Cf0024Col}}\\
     \end{center}
     \caption{The damping of a discharge wave (a) $C_f=0$, $\alpha=\infty$, (b) $C_f=0.000022$, $\alpha=15.15$, 
     (c) $C_f=0.000202$, $\alpha=5$ and (d) $C_f=0.005053$, $\alpha=1$. 
     The friction has been treated with either the semi-implicit method (SI) or the apparent topography method (AT)
 and $J=800$.}
     \label{fig-wave-damping}
\end{figure}

The following numerical values allow to plot figure \ref{fig-wave-damping-Cf0Col} to \ref{fig-wave-damping-Cf0024Col}.
For $C_f$, we took 4 different values ($C_f=0$, $C_f=0.000022$, $C_f=0.000202$ and $C_f=0.005053$,
corresponding to Womersley parameters   $\alpha=\infty$, $\alpha=15.15$, $\alpha=5$ and $\alpha=1$), 
 $J=50$, $100$, $200$, $400$, $800$ cells, $k=1.10^8 \text{Pa/m}$,
 $\rho=1060\text{m}^3$, $R_0=4.10^{-3}\text{m}$, $L=3 \text{m}$, $T_{end}=25. \text{s}$.
As initial conditions, we take a fluid at rest $
 Q(x,0)=0\text{m}^3/\text{s}$ 
and as input boundary condition
\begin{equation*}
Q_b(t)=Q_{amp}\sin(\omega t)\text{m}^3/\text{s}.
\end{equation*}
with $\omega=2\pi/T_{pulse}=2\pi/0.5\text{s}$ and $Q_{amp}=3.45.10^{-7}\text{m}^3/\text{s}$. As the flow is
 subcritical,the discharge is imposed at the outflow boundary thanks to \eqref{dampingDischarge} with $x=L$.

In this part, we insist on the comparison between first order and second order. So, we compare first order (HHL1) and second
 order scheme with both MUSCL (HLL MUSCL) and ENO (HLL ENO) linear reconstruction.
 When $C_f$ is not null, semi-implicit (SI) and apparent topography (AT) are compared. We should remark that as the friction
 increases the structure of the system \eqref{system1} changes and goes from a transport/ wave behavior to a diffusive behavior.
 The diffusive behavior comes from the fact that as $C_f$ increases, finally 
the momentum equation contains at leading order only the friction term and the elastic one:
$C_f Q_1/A_0 \sim -(A_0k/\rho) \partial_x R$ and so, using the conservation of mass:
$$\frac{\partial R}{\partial t} = D \frac{\partial^2 R}{\partial x^2}\mbox{  with  } 
D = \frac{k \pi R_0^2 }{2 \rho C_f},$$ 
which is consistent with \eqref{disp} at small $\omega$. However solving a "heat equation" with a
 "wave technique" seems to be adapted. Indeed, we notice a fair agreement between
 the numerics and the analytical solution in Figures \ref{fig-wave-damping-Cf0Col} to
 \ref{fig-wave-damping-Cf0024Col}. "Errors made" between numerical results and the exact solution of the linearized system are
 given for information with CPU time in tables \ref{tbl:TabErDamping2} and \ref{tbl:TabErDamping1}.

\subsection{Discussion on numerical precision}
\begin{figure}[htbp]
\begin{center}
       \includegraphics[angle=-90,width=0.48\textwidth]{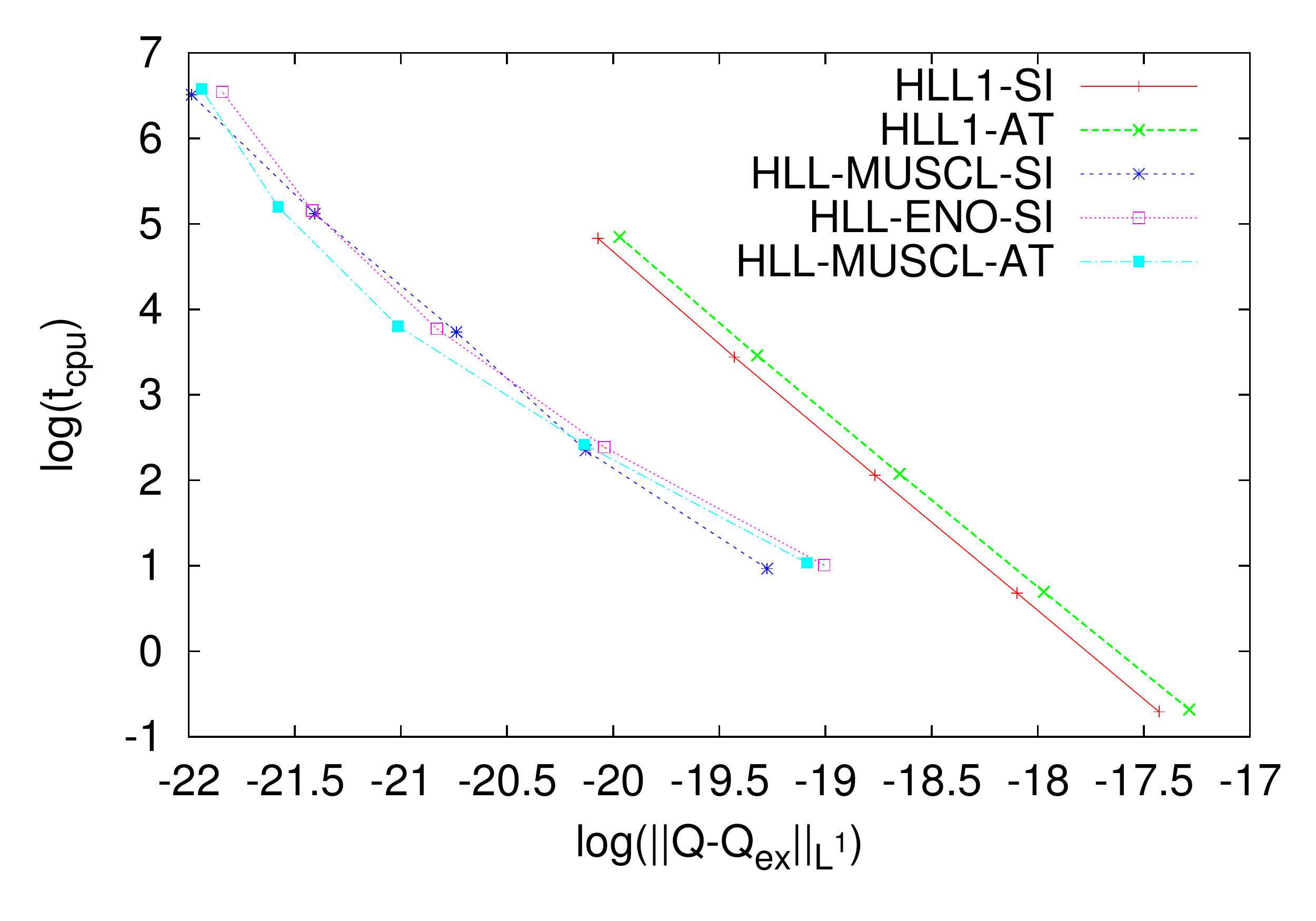}
     \caption{The wave "damping": The efficiency curves for $C_f=0.000202$.}
     \label{fig-efficiencydamping}
\end{center}
\end{figure}

Indeed, the notion of scheme order fully depends on the regularity of the solution. 
Moreover, calculating the error on the solution of the linearized system does not provide information from a numerical point of view
(sections \ref{we} and \ref{wd}). Our purpose is not to verify the second
 order accuracy of the numerical method (it has already been verified for the Shallow Water applications in
 \cite{Bouchut04,Audusse05,Delestre10}). But these errors computations allow us to verify that for these configurations
 the nonlinear system has a linear behavior.
 And this linear behavior is perfectly captured by the scheme and the physical model. Precision is increased by increasing the order
 of the numerical method (as shown on Figure \ref{fig-waveError} and as illustrated in tables \ref{tbl:TabErDamping2}
 and \ref{tbl:TabErDamping1}). The theoretical order of the schemes is not recovered for the wave test in section \ref{we}
 (0.78 instead of 1 and 1.15 instead of 2) this is due to the regularity of the solution of the initial system. Indeed, this solution
 is only continuous.
 In addition, the system has no term of order 2, then the interest in using a second order method
 is to save computing time for a given accuracy (see Figure \ref{fig-efficiencydamping}). The semi-implicit treatment and
 the apparent topography give closed results. However, when friction increases ($C_f=0.005053$), the semi-implicit method is
 closer to the linearised solution.

\begin{table}[htbp]
\begin{center}

{\small\begin{tabular}{|c|c|c|c|c|c|c|}  \cline{1-7}
 \multicolumn{1}{|c|}{$C_f=0$} &
\multicolumn{2}{c|}{HLL1} &
\multicolumn{2}{c|}{HLL MUSCL} &
\multicolumn{2}{c|}{HLL ENO} \\ \hline
\multicolumn{1}{|c|}{J} &
\(||Q-Q_{ex}||_{L^1}\) & \(t_{cpu}\) [s] &
\(||Q-Q_{ex}||_{L^1}\) & \(t_{cpu}\) [s] &
\(||Q-Q_{ex}||_{L^1}\) & \(t_{cpu}\) [s]\\ \hline
50  & 0.314E-7               & 0.49   & 0.104E-7      & 2.63          & 0.11E-7  & 2.7     \\ \hline
100  & 0.159E-7              & 1.93   & 0.508E-8      & 10.49         & 0.546E-8 & 10.7    \\ \hline
200  & 0.801E-8              & 7.7    & 0.238E-8      & 40.93         & 0.263E-8 & 42.6    \\ \hline
400  & 0.397E-8              & 30.93  & 0.103E-8      & 163.47        & 0.118E-8 & 170.32  \\ \hline
800 &  0.192E-8              & 123.28 & 0.36E-9       & 654.64        & 0.433E-9 & 681.26  \\ \hline
Regression & y=-1.007x-13.32 &        & y=-1.2x-13.59 &               & y=-1.155x-13.72 &      \\ \hline
\end{tabular}}
\end{center}

\caption{The wave "damping": $L^1$ errors on the discharge $Q$ and CPU times $t_{cpu}$ for $C_f=0$.}
\label{tbl:TabErDamping1}
\end{table}

 \section{Conclusion}
 \label{conclusion}
 
 In this paper we have considered the classical 1-D model of flow in an artery, we have presented a new numerical scheme and
 some numerical tests. The new proposed scheme has been obtained in following recent advances in the shallow-water Saint-Venant
 community. This community has been confronted to spurious flows induced by the change of topography in not well written
 schemes. In blood-flows, it corresponds to the treatment of the terms due to a variating initial shape of the artery, which
 arises is stenosis, aneurisms, or taper.
 To write the method, we have insisted on the fact that the effective conservative variables are: the artery area $A$ and the
 discharge $Q$, which was not clearly observed up to now. 
The obtained conservative system has been then discretized in order to use the property of equilibrium of "the man at eternal
 rest" analogous to the "lake at rest" in Saint-Venant equations. 
If this property is not preserved in the numerical stencil, spurious currents arise in the case of a variating vessel. For
 sake of illustration, if the terms are treated in a too simple way, we have exhibited a case with an aneurism which induces
 a non zero flux of blood $Q$. In a pulsatile case, the same configuration creates extra waves as well. An analogous of dam
 break (here a tourniquet) was used to validate the other terms of the dicretisation. Other less demanding examples have been
 performed: linear waves without and with damping in straight tubes. Good behaviour is obtained in those pertinent test cases
 that explore all the parts of the equations.

The well-balanced finite volume scheme that we propose preserves the volume of blood and avoids non physical case behavior.
 Thus we get a stable scheme and the accuracy is improved thanks to a second order reconstruction.
 The numerical method has been tested on examples which are more "gedanken" than physiological. 
This means that we only did a first necessary step. The next one is now to test more complex cases such as bifurcations,
 special boundary conditions and so on, in order to confront and fit to real practical clinical cases one of them being the
 case of aneurisms, in which as well we can introduce a source term in the mass equation in the case of leaks, we conclude
 that to look at those difficult problems, a careful treatment of the numerical scheme is important.\\

We thank ENDOCOM ANR for financial support, we thank Marie-Odile Bristeau and Christophe Berthon for fruitful discussions.\\

\bibliography{biblio}

\begin{thebibliography}{10}

\bibitem{Alastruey08}
J.~Alastruey, K.~H. Parker, J.~Peir\'o, and S.~J. Sherwin.
\newblock Lumped parameter outflow models for 1-d blood flow simulations:
  Effect on pulse waves and parameter estimation.
\newblock {\em Commun. Comp. Phys.}, 4:317--336, 2008.

\bibitem{Audusse99}
E.~Audusse.
\newblock Sch\'emas cin\'etiques pour le syst\`eme de {S}aint-{V}enant.
\newblock Master's thesis, Universit{\'e} Paris 6, 1999.

\bibitem{Audusse04c}
E.~Audusse, F.~Bouchut, M.-O. Bristeau, R.~Klein, and B.~Perthame.
\newblock A fast and stable well-balanced scheme with hydrostatic
  reconstruction for shallow water flows.
\newblock {\em SIAM J. Sci. Comput.}, 25(6):2050--2065, 2004.

\bibitem{Audusse05}
Emmanuel Audusse and Marie-Odile Bristeau.
\newblock A well-balanced positivity preserving {"}second-order{"} scheme for
  shallow water flows on unstructured meshes.
\newblock {\em Journal of Computational Physics}, 206:311--333, 2005.

\bibitem{Bermudez98}
Alfredo Berm\'udez, Alain Dervieux, Jean-Antoine Desideri, and M.~Elena
  V\'azquez.
\newblock Upwind schemes for the two-dimensional shallow water equations with
  variable depth using unstructured meshes.
\newblock {\em Computer Methods in Applied Mechanics and Engineering},
  155(1-2):49 -- 72, 1998.

\bibitem{Bermudez94}
Alfredo Bermudez and M.~Elena Vazquez.
\newblock Upwind methods for hyperbolic conservation laws with source terms.
\newblock {\em Computers \& Fluids}, 23(8):1049 -- 1071, 1994.

\bibitem{Berthon07}
Christophe Berthon and Fr\'ed\'eric Coquel.
\newblock Nonlinear projection methods for multi-entropies navier-stokes
  systems.
\newblock {\em Mathematics of Computation}, 76:1163--1194, 2007.

\bibitem{Bouchut04}
Fran\c{c}ois Bouchut.
\newblock {\em Nonlinear stability of finite volume methods for hyperbolic
  conservation laws, and well-balanced schemes for sources}, volume 2/2004.
\newblock Birkh{\"a}user Basel, 2004.

\bibitem{Bouchut09}
Fran\c{c}ois Bouchut and Tom\'as Morales.
\newblock A subsonic-well-balanced reconstruction scheme for shallow water
  flows.
\newblock http://www.math.ntnu.no/conservation/2009/032.html, 2009.
\newblock Preprint.

\bibitem{Bouchut07}
Fran\c{c}ois Bouchut, Haythem Ounaissa, and Beno\^it Perthame.
\newblock Upwinding of the source term at interfaces for euler equations with
  high friction.
\newblock {\em Computers and Mathematics with Applications}, 53:361--–375,
  2007.

\bibitem{Bristeau01}
Marie-Odile Bristeau and Beno\^it Coussin.
\newblock Boundary conditions for the shallow water equations solved by kinetic
  schemes.
\newblock Technical Report 4282, INRIA, October 2001.

\bibitem{Burguete08}
J.~Burguete, P.~Garc\'ia-Navarro, and J.~Murillo.
\newblock Friction term discretization and limitation to preserve stability and
  conservation in the 1d shallow-water model: Application to unsteady
  irrigation and river flow.
\newblock {\em International Journal for Numerical Methods in Fluids},
  58(4):403--425, 2008.

\bibitem{Castro08b}
Manuel~J. Castro, Philippe~G. LeFloch, Mar\'ia Luz Mu\~noz Ruiz, and Carlos
  Par\'es.
\newblock Why many theories of shock waves are necessary: Convergence error in
  formally path-consistent schemes.
\newblock {\em Journal of Computational Physics}, 227(17):8107--8129, 2008.

\bibitem{Castro07}
Manuel~J. Castro, Alberto Pardo~Milan\'es, and Carlos Par\'es.
\newblock Well-balanced numerical schemes based on a generalized hydrostatic
  reconstruction technique.
\newblock {\em Mathematical Models and Methods in Applied Sciences},
  17(12):2065--2113, 2007.

\bibitem{Cavallini08}
N.~Cavallini, V.~Caleffi, and V.~Coscia.
\newblock Finite volume and weno scheme in one-dimensional vascular system
  modelling.
\newblock {\em Computers and Mathematics with Applications}, 56:2382--2397,
  2008.

\bibitem{Cavallini10}
Nicola Cavallini and Vincenzo Coscia.
\newblock One-dimensional modelling of venous pathologies: Finite volume and
  weno schemes.
\newblock In Rolf Rannacher and Ad\'elia Sequeira, editors, {\em Advances in
  Mathematical Fluid Mechanics}, pages 147--170. Springer Berlin Heidelberg,
  2010.

\bibitem{Chinnayya99}
Ashwin Chinnayya and A.-Y. Le{R}oux.
\newblock A new {G}eneral {R}iemann {S}olver for the {S}hallow {W}ater
  equations, with friction and topography.
\newblock http://www.math.ntnu.no/conservation/1999/021.html, 1999.
\newblock Preprint.

\bibitem{saintvenant71}
{A}dh\'emar {J}ean-{C}laude de~{S}aint {V}enant.
\newblock Th\'eorie du mouvement non-permanent des eaux, avec application aux
  crues des rivi\`eres et \`a l'introduction des mar\'ees dans leur lit.
\newblock {\em Comptes Rendus de l'Acad\'emie des {S}ciences}, 73:147--154,
  1871.

\bibitem{Delestre10b}
Olivier Delestre.
\newblock {\em Simulation du ruissellement d'eau de pluie sur des surfaces
  agricoles/ rain water overland flow on agricultural fields simulation}.
\newblock PhD thesis, Universit\'e d'Orl\'eans (in French), available from TEL:
  tel.archives-ouvertes.fr/INSMI/tel-00531377/fr, July 2010.

\bibitem{Delestre09b}
Olivier Delestre, St\'ephane Cordier, Fran\c{c}ois James, and Fr\'ed\'eric
  Darboux.
\newblock Simulation of rain-water overland-flow.
\newblock In {\em Proceedings of the 12th International Conference on
  Hyperbolic Problems, \rm University of Maryland, College Park (USA), 2008, E.
  Tadmor, J.-G. Liu and A. Tzavaras Eds., Proceedings of Symposia in Applied
  Mathematics 67, Amer. Math. Soc., 537--546}, 2009.

\bibitem{Delestre09}
Olivier Delestre and Fran\c{c}ois James.
\newblock Simulation of rainfall events and overland flow.
\newblock In {\em Proceedings of X International Conference Zaragoza-Pau on
  Applied Mathematics and Statistics, \rm Jaca, Spain, september 2008,
  Monograf\'ias Matem\'aticas Garc\'ia de Galdeano}, 2009.

\bibitem{Delestre10}
Olivier Delestre and Fabien Marche.
\newblock A numerical scheme for a viscous shallow water model with friction.
\newblock {\em J. Sci. Comput.}, DOI 10.1007/s10915-010-9393-y, 2010.

\bibitem{Dutykh07}
Denys Dutykh.
\newblock {\em Mod\'elisation math\'ematique des tsunamis}.
\newblock PhD thesis, \'Ecole normale sup\'erieure de Cachan, 2007.

\bibitem{Fernandez07}
M.~A. Fern\'andez, J.-F. Gerbeau, and C.~Grandmont.
\newblock A projection semi-implicit scheme for the coupling of an elastic
  structure with an incompressible fluid.
\newblock {\em International Journal for Numerical Methods in Engineering},
  69(4):794--821, 2007.

\bibitem{Fernandez04b}
Miguel~\'Angel Fern\'andez, Vuk Mili{\v s}i\'c, and Alfio Quarteroni.
\newblock Analysis of a geometrical multiscale blood flow model based on the
  coupling of ode's and hyperbolic pde's.
\newblock Technical Report 5127, INRIA, February 2004.

\bibitem{Fiedler00}
R.~F. Fiedler and J.~A. Ramirez.
\newblock A numerical method for simulating discontinuous shallow flow over an
  infiltrating surface.
\newblock {\em International Journal for Numerical Methods in Fluids},
  32:219--240, 2000.

\bibitem{Formaggia06}
Luca Formaggia, Daniele Lamponi, Massimiliano Tuveri, and Alessandro Veneziani.
\newblock Numerical modeling of 1d arterial networks coupled with a lumped
  parameters description of the heart.
\newblock {\em Computer Methods in Biomechanics and Biomedical Engineering},
  9:273--288, 2006.

\bibitem{Fullana09}
Jose-Maria Fullana and St\'ephane Zaleski.
\newblock A branched one-dimensional model of vessel networks.
\newblock {\em J. Fluid. Mech.}, 621:183--204, 2009.

\bibitem{Gallouet03}
T.~Gallou\"et, J.-M. H\'erard, and N.~Seguin.
\newblock Some approximate godunov schemes to compute shallow-water equations
  with topography.
\newblock {\em Computers \& Fluids}, 32:479--513, 2003.

\bibitem{Gosse00}
L.~Gosse.
\newblock A well-balanced flux-vector splitting scheme designed for hyperbolic
  systems of conservation laws with source terms.
\newblock {\em Computers \& Mathematics with Applications}, 39(9-10):135 --
  159, 2000.

\bibitem{Goutal02}
N.~Goutal and F.~Maurel.
\newblock A finite volume solver for 1{D} shallow-water equations applied to an
  actual river.
\newblock {\em International Journal for Numerical Methods in Fluids},
  38:1--19, 2002.

\bibitem{Greenberg96}
J.~M. Greenberg and A.-Y. Le{R}oux.
\newblock A well-balanced scheme for the numerical processing of source terms
  in hyperbolic equation.
\newblock {\em SIAM Journal on Numerical Analysis}, 33:1--16, 1996.

\bibitem{Harten87}
Ami Harten, Bjorn Engquist, Stanley Osher, and Sukumar~R. Chakravarthy.
\newblock Uniformly {H}igh {O}rder {A}ccurate {E}ssentially {N}on-oscillatory
  {S}chemes, {III}.
\newblock {\em Journal of Computational Physics}, 71:231--303, August 1987.

\bibitem{Harten86}
Ami Harten, Stanley Osher, Bjorn Engquist, and Sukumar~R. Chakravarthy.
\newblock Some results on uniformly high-order accurate essentially
  nonoscillatory schemes.
\newblock {\em Applied Numerical Mathematics}, 2(3-5):347 -- 377, 1986.
\newblock Special Issue in Honor of Milt Rose's Sixtieth Birthday.

\bibitem{Harten83}
Amiram Harten, Peter~D. Lax, and Bram van Leer.
\newblock On upstream differencing and godunov-type schemes for hyperbolic
  conservation laws.
\newblock {\em SIAM Review}, 25(1):35--61, January 1983.

\bibitem{Hou94}
Thomas~Y. Hou and Philippe~G. LeFloch.
\newblock Why nonconservative schemes converge to wrong solutions: error
  analysis.
\newblock {\em Mathematics of Computation}, 62(206):497--530, April 1994.

\bibitem{Jin01a}
Shi Jin.
\newblock A steady-state capturing method for hyperbolic systems with
  geometrical source terms.
\newblock {\em M2AN}, 35(4):631--645, July 2001.

\bibitem{Kirkman03}
Robert Kirkman, Tony Moore, and Charlie Adlard.
\newblock {\em The Walking Dead}.
\newblock Image Comics, 2003.

\bibitem{Kundu04}
Pijush~K. Kundu and Ira~M. Cohen.
\newblock {\em Fluid Mechanics, Third Edition}.
\newblock 2004.

\bibitem{Kundu08}
Pijush~K. Kundu and Ira~M. Cohen.
\newblock {\em Fluid Mechanics, Fourth Edition}.
\newblock 2008.

\bibitem{Kurganov02}
Alexander Kurganov and Doron Levy.
\newblock Central-upwind schemes for the saint-venant system.
\newblock {\em Mathematical Modelling and Numerical Analysis}, 36:397--425,
  2002.

\bibitem{Lagree00}
Pierre-Yves Lagr\'ee.
\newblock An inverse technique to deduce the elasticity of a large artery.
\newblock {\em The European Physical Journal}, 9:153--163, 2000.

\bibitem{Lagree05}
Pierre-Yves Lagr\'ee and Sylvie Lorthois.
\newblock The rns/prandtl equations and their link with other asymptotic
  descriptions. application to the computation of the maximum value of the wall
  shear stress in a pipe.
\newblock {\em Int. J. Eng. Sci.}, 43/3--4:352--378, 2005.

\bibitem{Lagree96}
Pierre-Yves Lagr\'ee and Maurice Rossi.
\newblock Etude de l'\'ecoulement du sang dans les art\`eres: effets
  nonlin\'eaires et dissipatifs.
\newblock {\em C. R. Acad. Sci. Paris, t322, S\'erie II b}, pages 401--408,
  1996.

\bibitem{Lee10}
Sang-Heon Lee and Nigel~G. Wright.
\newblock Simple and efficient solution of the shallow water equations with
  source terms.
\newblock {\em International Journal for Numerical Methods in Fluids},
  63:313--340, 2010.

\bibitem{LeVeque92}
Randall~J. LeVeque.
\newblock {\em Numerical methods for conservation laws}.
\newblock Lectures in mathematics ETH Zurich. 1992.

\bibitem{LeVeque98}
Randall~J. LeVeque.
\newblock Balancing source terms and flux gradients in high-resolution godunov
  methods: The quasi-steady wave-propagation algorithm.
\newblock {\em Journal of Computational Physics}, 146(1):346 -- 365, 1998.

\bibitem{LeVeque02}
Randall~J. LeVeque.
\newblock {\em Finite volume methods for hyperbolic problems}.
\newblock Cambridge Texts in Applied Mathematics. Cambridge University Press,
  Cambridge, 2002.

\bibitem{Liang09b}
Qiuhua Liang and Fabien Marche.
\newblock Numerical resolution of well-balanced shallow water equations with
  complex source terms.
\newblock {\em Advances in Water Resources}, 32(6):873 -- 884, 2009.

\bibitem{Lighthill78}
James Lighthill.
\newblock {\em Waves in Fluids}.
\newblock 1978.

\bibitem{LukacovaMedvidova06}
M.~Luk\'a{\v c}ov\'a-Medvid'ov\'a and U.~Teschke.
\newblock Comparison study of some finite volume and finite element methods for
  the shallow water equations with bottom topography and friction terms.
\newblock {\em J. Appl. Mech. Math. (ZAMM)}, 86(11):874--891, 2006.

\bibitem{LukacovaMedvidova05b}
M.~Luk\'a{\v c}ov\'a-Medvid'ov\'a and Z.~Vlk.
\newblock Well-balanced finite volume evolution {G}alerkin methods for the
  shallow water equations with source terms.
\newblock {\em International Journal for Numerical Methods in Fluids},
  47(10-11):1165--1171, 2005.

\bibitem{Mangeney07}
A.~Mangeney, F.~Bouchut, N.~Thomas, J.~P. Vilotte, and M.-O. Bristeau.
\newblock Numerical modeling of self-channeling granular flow and of their
  levee-channel deposits.
\newblock {\em Journal of Geophysical Research}, 112(F2):F02017, May 2007.

\bibitem{MangeneyCastelnau05}
A.~Mangeney-Castelnau, F.~Bouchut, J.~P. Vilotte, E.~Lajeunesse, A.~Aubertin,
  and M.~Pirulli.
\newblock On the use of saint-venant equations to simulate the spreading of a
  granular mass.
\newblock {\em Journal of Geophysical Research}, 110:B09103, 2005.

\bibitem{Marchandise08}
Emilie Marchandise, Nicolas Chevaugeon, and Jean-Fran\c{c}ois Remacle.
\newblock Spatial and spectral superconvergence of discontinuous galerkin
  method for hyperbolic problems.
\newblock {\em Journal of Computational and Applied Mathematics},
  215(2):484--494, 2008.
\newblock Proceedings of the Third International Conference on Advanced
  Computational Methods in Engineering (ACOMEN 2005)., Proceedings of the Third
  International Conference on Advanced Computational Methods in Engineering
  (ACOMEN 2005).

\bibitem{Marche07b}
F.~Marche, P.~Bonneton, P.~Fabrie, and N.~Seguin.
\newblock Evaluation of well-balanced bore-capturing schemes for 2{D} wetting
  and drying processes.
\newblock {\em International Journal for Numerical Methods in Fluids},
  53:867--894, 2007.

\bibitem{Marche08}
Fabien Marche and Christophe Berthon.
\newblock A positive preserving high order vfroe scheme for shallow water
  equations: A class of relaxation schemes.
\newblock {\em {SIAM} {J}ournal on {S}cientific {C}omputing}, 30(5):pp.
  2587--2612, August 2008.
\newblock 75{M}12, 35{L}65, 65{M}12.

\bibitem{Martin05}
Vincent Martin, Fran\c{c}ois Cl\'ement, Astrid Decoene, and Jean-Fr\'ed\'eric
  Gerbeau.
\newblock Parameter identification for a one-dimensional blood flow model.
\newblock In Eric Canc\`es \& Jean-Fr\'ed\'eric Gerbeau, editor, {\em ESAIM:
  PROCEEDINGS}, volume~14, pages 174--200, September 2005.

\bibitem{Noelle06}
Sebastian Noelle, Normann Pankratz, Gabriella Puppo, and Jostein~R. Natvig.
\newblock Well-balanced finite volume schemes of arbitrary order of accuracy
  for shallow water flows.
\newblock {\em Journal of Computational Physics}, 213(2):474 -- 499, 2006.

\bibitem{Noelle07}
Sebastian Noelle, Yulong Xing, and Chi~Wang Shu.
\newblock High-order well-balanced finite volume weno schemes for shallow water
  equation with moving water.
\newblock {\em Journal of Computational Physics}, 226(1):29 -- 58, 2007.

\bibitem{Ockendon83}
Hilary Ockendon and Alan~B. Tayler.
\newblock {\em Inviscid Fluid Flows}.
\newblock 1983.

\bibitem{Olufsen99}
Mette~S. Olufsen.
\newblock Structured tree outflow condition for blood flow in larger systemic
  arteries.
\newblock {\em American Journal of Physiology -- Heart and Circulatory
  Physiology}, 276:257--268, 1999.

\bibitem{Olufsen00}
Mette~S. Olufsen, Charles~S. Peskin, Won~Yong Kim, Erik~M. Pedersen, Ali Nadim,
  and Jesper Larsen.
\newblock Numerical simulation and experimental validation of blood flow in
  arteries with structured-tree outflow conditions.
\newblock {\em Annals of Biomedical Engineering}, 28:1281--1299, 2000.

\bibitem{Parker09}
Kim~H. Parker.
\newblock A brief history of arterial wave mechanics.
\newblock {\em Med Biol Eng Comput}, 47:111--118, February 2009.

\bibitem{Pedley80}
T.~J. Pedley.
\newblock {\em The Fluid Mechanics of Large Blood Vessel}.
\newblock 1980.

\bibitem{Perthame01}
B.~Perthame and C.~Simeoni.
\newblock A kinetic scheme for the {S}aint-{V}enant system with a source term.
\newblock {\em Calcolo}, 38:201--231, 2001.
\newblock 10.1007/s10092-001-8181-3.

\bibitem{Pindera09}
Maciej~Z. Pindera, Hui Ding, Mahesh~M. Athavale, and Zhijian Chen.
\newblock Accuracy of 1d microvascular flow models in the limit of low reynolds
  numbers.
\newblock {\em Microvascular Research}, 77(3):273--280, 2009.

\bibitem{Popinet11}
St\'ephane Popinet.
\newblock Quadtree-adaptative tsunami modelling.
\newblock {\em Ocean Dynamics}, pages 1--25, May 2011.

\bibitem{Saito11}
Masashi Saito, Yuki Ikenaga, Mami Matsukawa, Yoshiaki Watanabe, Takaaki Asada,
  and Pierre-Yves Lagr\'ee.
\newblock One-dimensional model for propagation of a pressure wave in a model
  of the human arterial network: Comparison of theoretical and experimental.
\newblock {\em Journal of Biomechanical Engineering}, 133, December 2011.

\bibitem{Sampson05}
Joe Sampson, Alan Easton, and Manmohan Singh.
\newblock Modelling the effect of proposed channel deepening on the tides in
  {P}ort {P}hillip {B}ay.
\newblock 2005.

\bibitem{Sherwin03}
S.~J. Sherwin, L.~Formaggia, J.~Peir\'o, and V.~Franke.
\newblock Computational modelling of 1d blood flow with variable mechanical
  properties and its application to the simulation of wave propagation in the
  human arterial system.
\newblock {\em International Journal for Numerical Methods in Fluids},
  43:673--700, 2003.

\bibitem{Shu88}
C.-W. Shu and Stanley Osher.
\newblock Efficient implementation of essentially non-oscillatory
  shock-capturing schemes.
\newblock {\em Journal of Computational Physics}, 77(2):439 -- 471, 1988.

\bibitem{Stergiopulos92}
N.~Stergiopulos, D.~F. Young, and T.~R. Rogge.
\newblock Computer simulation of arterial flow with applications to arterial
  and aortic stenoses.
\newblock {\em J. Biomechanics}, 25(12):1477--1488, 1992.

\bibitem{Stoker57}
J.~J. Stoker.
\newblock {\em Water Waves: The Mathematical Theory with Applications}.
\newblock Interscience Publishers, New York, USA, 1957.

\bibitem{Tatard08}
L.~Tatard, O.~Planchon, J.~Wainwright, G.~Nord, D.~Favis-Mortlock, N.~Silvera,
  O.~Ribolzi, M.~Esteves, and Chi-Hua Huang.
\newblock Measurement and modelling of high-resolution flow-velocity data under
  simulated rainfall on a low-slope sandy soil.
\newblock {\em Journal of Hydrology}, 348(1-2):1--12, January 2008.

\bibitem{Thanh08}
M.~D. Thanh, Md. {F}azlul {K}arim, and A.~I.~Md. Ismail.
\newblock Well-balanced scheme for shallow water equations with arbitrary
  topography.
\newblock {\em Int. J. Dynamical Systems and Differential Equations},
  1(3):196--204, 2008.

\bibitem{Valiani99a}
A.~Valiani, V.~Caleffi, and A.~Zanni.
\newblock Finite volume scheme for {2D} shallow-water equations. {A}pplication
  to {M}alpasset dam-break.
\newblock In {\em the $4^{th}$ {CADAM} {W}orkshop, {Z}aragoza}, pages 63--94,
  1999.

\bibitem{Valiani02}
A.~Valiani, V.~Caleffi, and A.~Zanni.
\newblock Case {S}tudy : {M}alpasset {D}am-{B}reak {S}imulation using a
  {T}wo-{D}imensional {F}inite {V}olume {M}ethods.
\newblock {\em Journal of {H}ydraulic {E}ngineering}, 128(5):460--472, May
  2002.

\bibitem{vandeVosse03b}
F.N. van~de Vosse, J.~de~Hart, C.H.G.A. van Oijen, and D.~Bessems.
\newblock Finite-element-based computational methods for cardiovascular
  fluid-structure interaction.
\newblock {\em Journal of Engineering Mathematics}, 47:335--368, 2003.

\bibitem{vanLeer79}
Bram van {L}eer.
\newblock {T}owards the ultimate conservative difference scheme. {V}. {A}
  second-order sequel to {G}odunov's method.
\newblock {\em Journal of Computational Physics}, 32(1):101 -- 136, 1979.

\bibitem{VanSteenhoven86}
A.~A. Van~Steenhoven and M.~E.~H. Van~Dongen.
\newblock Model studies of the aortic pressure rise just after valve closure.
\newblock {\em J. Fluid. Mech.}, 166:93--113, 1986.

\bibitem{Vazquez-Cendon99}
M.~E. V\`azquez-Cend\`on.
\newblock Improved treatment of source terms in upwind schemes for the shallow
  water equations in channels with irregular geometry.
\newblock {\em Journal of Computational Physics}, 148:497--526, 1999.

\bibitem{Wibmer04}
Michael Wibmer.
\newblock {\em One-dimensional Simulation of Arterial Blood Flow with
  Applications}.
\newblock PhD thesis, eingereicht an der Technischen Universitat Wien --
  Fakultat fur Technische Naturwissenschaften und Informatik, January 2004.

\bibitem{Willemet11}
Marie Willemet, Val\'erie Lacroix, and Emilie Marchandise.
\newblock Inlet boundary conditions for blood flow simulations in truncated
  arterial networks.
\newblock {\em Journal of Biomechanics}, 44(5):897--903, 2011.

\bibitem{Womersley55}
J.R. Womersley.
\newblock On the oscillatory motion of a viscous liquid in thin-walled elastic
  tube: I.
\newblock {\em Philos. Mag.}, 46:199--221, 1955.

\bibitem{Xiu07}
Dongbin Xiu and Spencer~J. Sherwin.
\newblock Parametric uncertainty analysis of pulse wave propagation in a model
  of a human arterial network.
\newblock {\em Journal of Computational Physics}, 226:1385--1407, 2007.

\bibitem{Zagzoule91}
M.~Zagzoule, J.~Khalid-Naciri, and J.~Mauss.
\newblock Unsteady wall shear stress in a distensible tube.
\newblock {\em J. Biomechanics}, 24(6):435--439, 1991.

\bibitem{Zagzoule86}
Mokhtar Zagzoule and Jean-Pierre Marc-Vergnes.
\newblock A global mathematical model of the cerebral circulation in man.
\newblock {\em J. Biomechanics}, 19(12):1015--1022, 1986.

\end{thebibliography}
\bibliographystyle{plain}

\end{document}